\documentclass[screen,authorversion,nonacm]{acmart}

\AtBeginDocument{%
  }

\setcopyright{acmlicensed}
\copyrightyear{2024}
\acmYear{2024}
\acmDOI{XXXXXXX.XXXXXXX}  %

\acmJournal{TOMS}
\acmVolume{0}
\acmNumber{0}
\acmArticle{0}
\acmMonth{0}

\usepackage{lipsum}
\usepackage{svg}
\usepackage{mathtools}   %
\usepackage[10pt]{moresize}
\usepackage[htt]{hyphenat}
\usepackage[utf8]{inputenc}

\usepackage{adjustbox}
\usepackage{multirow}    %
\usepackage{pifont}      %

\usepackage{cleveref}    %
\usepackage{makecell}
\usepackage{soul}

\ifpdf
  \DeclareGraphicsExtensions{.eps,.pdf,.png,.jpg}
\else
  \DeclareGraphicsExtensions{.eps}
\fi

\usepackage{julia-mono-listings}
\usepackage[frozencache,cachedir=.]{minted}
\usepackage{tcolorbox}
\usepackage{etoolbox}

\ifdefined\directlua
  \usepackage{fontspec}
\else
  \usepackage[T1]{fontenc}
  \usepackage[nomath]{lmodern}
\fi

\usepackage{tikz}
\usetikzlibrary{arrows.meta, positioning, shapes, fit, backgrounds}

\tikzset{
    box/.style={rectangle, draw=black, thick, text centered, minimum height=1.2em, inner sep=4pt},
    bluebox/.style={box, draw=cyan, dashed, fill=cyan!10},
    graybox/.style={box, draw=gray, dashed, fill=gray!10},
    purplebox/.style={box, draw=violet, dashed, fill=violet!10},
    greenbox/.style={box, draw=green!60!black, dashed, fill=green!10},
    orangebox/.style={box, draw=orange, dashed, fill=orange!10},
    arrow/.style={-Stealth, thick},
}

\newcommand{\func}[2]{#1\left(#2\right)}

\newcommand{\nlsolve}{\texttt{NonlinearSolve.jl}}

\renewcommand{\descriptionlabel}[1]{\hspace*{3em}\ul{#1}}

\newcommand{\cmark}{\textcolor{green!70!black}{\ding{51}}} %
\newcommand{\xmark}{\textcolor{red}{\ding{55}}}            %

\DeclarePairedDelimiter{\norm}{\lVert}{\rVert}

\DeclareMathOperator*{\argmin}{arg\,min}


\begin{document}

\setlength{\footskip}{40pt}
\makeatletter
\fancyfoot{} %
\fancyfoot[RE,LO]{\textit{Distribution Statement A: Approved for public release. Distribution is Unlimited.}}
\renewcommand{\@oddfoot}{\if@ACM@sigchiamode\hfill\thepage\hfill\else\hfill\thepage\hfill\fi}
\renewcommand{\@evenfoot}{\if@ACM@sigchiamode\hfill\thepage\hfill\else\hfill\thepage\hfill\fi}
\pagestyle{fancy}
\makeatother

\title[High-Performance Solvers for Systems of Nonlinear Equations]{NonlinearSolve.jl: High-Performance and Robust Solvers for Systems of Nonlinear Equations in Julia}

\author{Avik Pal}
\email{avikpal@mit.edu}
\orcid{0002-3938-7375}

\author{Flemming Holtorf}
\affiliation{%
  \institution{Massachusetts Institute of Technology}
  \city{Cambridge}
  \state{MA}
  \country{USA}
}

\author{Axel Larsson}
\affiliation{%
  \institution{Princeton University}
  \city{Princeton}
  \state{NJ}
  \country{USA}
}

\author{Torkel Loman}
\author{Utkarsh}
\author{Frank Sch\"afer}
\affiliation{%
  \institution{Massachusetts Institute of Technology}
  \city{Cambridge}
  \state{MA}
  \country{USA}
}

\author{Qingyu Qu}
\affiliation{%
    \institution{Zhejiang University}
    \country{China}
}

\author{Alan Edelman}
\email{edelman@mit.edu}

\author{Chris Rackauckas}
\email{crackauc@mit.edu}
\affiliation{%
  \institution{Massachusetts Institute of Technology}
  \city{Cambridge}
  \state{MA}
  \country{USA}
}

\renewcommand{\shortauthors}{Pal, Holtorf, Larsson, Loman, Sch\"after, Qu, Edelman, Rackauckas}
\acmArticleType{Research}
\acmCodeLink{https://github.com/SciML/NonlinearSolve.jl}
\acmDataLink{https://zenodo.org/records/10397607}

\authorsaddresses{Corresponding author: Avik Pal,
\href{mailto:avikpal@mit.edu}{avikpal@mit.edu};
Alan Edelman,
\href{mailto:edelman@mit.edu}{edelman@mit.edu};
Chris Rackauckas,
\href{mailto:crackauc@mit.edu}{crackauc@mit.edu};
CSAIL, MIT, Cambridge, MA, USA, 02139}

\begin{abstract}
    Efficiently solving nonlinear equations underpins numerous scientific and engineering disciplines, yet scaling these solutions for challenging system models remains a challenge. This paper presents \nlsolve{} -- a suite of high-performance open-source nonlinear equation solvers implemented natively in the Julia programming language. \nlsolve{} distinguishes itself by offering a unified API that accommodates a diverse range of solver specifications alongside features such as automatic algorithm selection based on runtime analysis, support for GPU-accelerated computation through static array kernels, and the utilization of sparse automatic differentiation and Jacobian-free Krylov methods for large-scale problem-solving. Through rigorous comparison with established tools such as \texttt{PETSc SNES}, \texttt{Sundials KINSOL}, and \texttt{MINPACK}, \nlsolve{} demonstrates robustness and efficiency, achieving significant advancements in solving nonlinear equations while being implemented in a high-level programming language. The capabilities of \nlsolve{} unlock new potentials in modeling and simulation across various domains, making it a valuable addition to the computational toolkit of researchers and practitioners alike.
\end{abstract}

\maketitle

\section{Introduction}\label{sec:introduction}

Nonlinear equations are ubiquitous across the sciences, governing the behavior of challenging systems in domains ranging from physics to biology to machine learning. Solving steady state (or fixed-point) equations, implicit ordinary differential equations (ODE)~\cite{wanner1996solving}, semi-explicit differential algebraic equations (DAEs), boundary value differential equations~\cite{Enright1996RungeKuttaSW, deuflhard1983multiple, ascher1995numerical, zhang2012improved}, and memory-efficient deep learning~\cite{baideep2019, baimultiscale2020, pal2023continuous, pal2023efficient} all boil down to solving nonlinear equations. Despite their prevalence in real-world problems, finding solutions to nonlinear equations efficiently remains a major challenge.

This paper presents \nlsolve{} - an open-source, high-performance nonlinear equation-solving framework implemented in the Julia programming language~\cite{bezanson2017julia}. \nlsolve{} provides a flexible, easy-to-use interface for specifying and solving general nonlinear equations. Under the hood, it employs various state-of-the-art techniques, such as automatic sparsity exploitation, fast automatic differentiation (AD), and Jacobian-free methods to solve large systems reliably. \nlsolve{} is part of the \href{https://sciml.ai/}{Scientific Machine Learning (SciML) package ecosystem of Julia} and is, in turn, extended by various packages to solve domain-specific problems. The key capabilities and contributions of \nlsolve{} include:
{\sloppy
\begin{itemize}
    \item \ul{Unified API for Rapid Experimentation with Nonlinear Solver Options}: Users can seamlessly switch between solver algorithms (including the same algorithm from different software packages), line search methods, trust-region schemes, automatic differentiation backends, linear solvers (including Krylov methods), and sparsity detection algorithms. This includes enabling sparsity and AD support for external packages like \texttt{PETSc} via the same API.

    \item \ul{Smart Polyalgorithmic Solver Selection for Robustness}: Automated runtime checks select appropriate internal solver methods and parameters, balancing speed and reliability.

    \item \ul{High-performance Algorithms for Small Problems}: Specialized non-allocating kernels solve small systems extremely efficiently -- outperforming existing software.

    \item \ul{Automatic Sparsity Exploitation for Large Systems}: Approximation techni\-ques automatically detect sparsity patterns in system Jacobians. Colored sparse matrix computations coupled with sparse linear solvers accelerate both derivative calculations and overall solve times.

    \item \ul{Jacobian-Free Krylov Methods}: Matrix-free linear solvers avoid explicit Jacobian materialization, reducing memory overhead and enabling the solution of otherwise computationally infeasible systems.

    \item \ul{Composable Blocks for Building Solvers}: All solvers in \texttt{NonlinearSolve.jl} are built using a composition of fundamental building blocks -- Descent Direction, Linear Solver, Jacobian Computation Algorithm, and Globalization Strategy.

    \item \ul{Symbolic Tearing of Nonlinear Systems}: \texttt{ModelingToolkit.jl}~\cite{ma2021modelingtoolkit} automatically performs tearing of nonlinear systems for \nlsolve{} and reduces the computational complexity of solving those systems.
\end{itemize}
}
We demonstrate the capabilities of \nlsolve{} with several numerical experiments. We perform numerical experiments and compare performance against \texttt{Sundials}, \texttt{MINPACK}, \texttt{PETSc}, and \texttt{NLsolve} in \Cref{sec:results} and demonstrate \nlsolve's ability to reliably solve root-finding problems. Furthermore, we show \nlsolve's applicability to large problems by benchmarking it on applications like initializing a DAE battery model [\Cref{subsec:battery_model}] and steady-state partial differential equation (PDE) solving [\Cref{subsec:brusselator}]. \nlsolve{} enables fast yet robust nonlinear equation solving on ill-conditioned problems. Its efficiency, flexibility, and seamless GPU support unlock new modeling and simulation capabilities across application domains.

\subsection{Julia as the Programming Language of Choice}\label{subsec:julia-prog-lang}

Julia~\cite{bezanson2017julia} is a programming language designed for technical computing, offering a unique combination of high-level syntax and high-performance execution. This makes it especially suited for numerical tasks such as nonlinear root-finding, where performance and composability are both critical. Julia’s ecosystem supports first-class automatic differentiation (AD), efficient linear algebra, and solver frameworks, enabling the rapid development of robust nonlinear solvers. The language's performance stems from its just-in-time (JIT) compilation strategy powered by LLVM~\cite{lattner2004llvm}. This allows numerical code—including AD routines and Jacobian computations—to run with performance comparable to C or Fortran, while maintaining a concise, high-level interface.

Multiple dispatch is a core feature of the language, allowing algorithm specialization based on types. This makes it easy to swap between dense, sparse, and matrix-free solvers, or to switch between forward- and reverse-mode AD backends, without changing user-facing code. AD frameworks like \texttt{ForwardDiff.jl}, \texttt{Zygote.jl}, and \texttt{Enzyme.jl} integrate natively with scientific computing packages, enabling accurate and efficient derivative computation even for large systems. GPU support via the JuliaGPU ecosystem~\cite{besard2018juliagpu,churavy2021juliagpu,besard2019rapid} further extends this flexibility, offering both high-level abstractions and low-level kernel control. This allows the same root-finding code to target CPUs, GPUs, or distributed environments with minimal changes. Furthermore, Julia's native support for calling C and Fortran, enables us to use state-of-the-art external solver packages like \texttt{Sundials} and \texttt{PETSc}. Julia's design principles—performance, composability, and expressiveness—make it an ideal foundation for building scalable and modular nonlinear root-finding tools with AD support.

\vspace{-1em}
\subsection{Comparison to Existing Software}\label{subsec:existing-software}

While there are several nonlinear solver suites, each with its advantages, they often have shortcomings in terms of flexibility, ease of use, performance, and scalability limitations when applied to real-world nonlinear problems. For example, \texttt{Sundials KINSOL}~\cite{gardner2022sundials, hindmarsh2005sundials} is a C library designed for solving nonlinear systems and supports various solver configurations through user-defined settings. However, its primary approach, which relies on a modified Newton specialized for implicit solvers, can be restrictive for certain applications, and techniques like Jacobian reuse often fail even for simple test problems [\Cref{subsec:23_test_problems}]. \texttt{MINPACK}~\cite{more1980user, cminpack} forms the foundation for widely used packages like \texttt{SciPy}~\cite{2020SciPy-NMeth} and \texttt{MATLAB fsolve}~\cite{OptimizationToolbox}; however, this lacks the support of tools like automatic differentiation, sparsity detection, and Krylov methods, is often more than $100$ times slower than our solvers [\Cref{fig:brusselator-scaling-benchmarks}] and fails for numerically challenging problems [\Cref{subsec:battery_model}]. \texttt{PETSc SNES}~\cite{balay2019petsc} is a specialized set of solvers for large systems in a distributed setup; however, they lack the support for sparse AD and instead can compute sparse Jacobians exclusively via colored finite differencing. \texttt{Optimistix}~\cite{rader2024optimistix} solves systems of nonlinear equations for machine learning-focused applications and, as such, lacks support for any form of sparsity handling that is uncommon in their target domain. We draw much inspiration from pre-existing Julia packages like \texttt{NLsolve.jl} and \texttt{NLSolvers.jl}~\cite{Optim.jl-2018}; however, they have limited flexibility and lack specialized linear solvers as evident from their slower relative performance [\Cref{sec:results}]. In contrast to existing software, \nlsolve{} is designed to be robust, performant, and scalable in general while having controls to be customized to specific applications.

\subsection{Domain Specific Extensions of \nlsolve}\label{subsec:nlsolve-extensions}

In addition to its core functionality for solving general nonlinear systems, \nlsolve{} is extended by several domain-specific packages that integrate nonlinear solvers into broader modeling workflows. For example, \texttt{SteadyStateDiffEq.jl} leverages \nlsolve{} to compute steady states of dynamical systems, with applications in systems biology and chemical reaction network modeling. Notably, \texttt{Catalyst.jl}~\cite{CatalystPLOSCompBio2023} uses this functionality to analyze steady states in biochemical systems such as dimerization reactions. In machine learning, \texttt{DeepEquilibriumNetworks.jl}~\cite{pal2023continuous} builds on \nlsolve{}’s sensitivity analysis capabilities to train deep equilibrium models, where inference corresponds to solving a nonlinear system. In numerical differential equations, \texttt{DiffEqCallbacks.jl} provides manifold-aware callbacks, such as \texttt{ManifoldCallback}, which employ \nlsolve{} to enforce algebraic constraints~\cite{hairer2006geometric} in problems like the Kepler system, preserving energy and angular momentum~\cite{cordani2003kepler}. Finally, \texttt{OrdinaryDiffEq.jl} and other solvers in the SciML ecosystem rely on \nlsolve{} for initializing and advancing implicit differential-algebraic equation (DAE) solvers, as demonstrated in the DAE battery model in \Cref{subsec:battery_model}. These integrations highlight \nlsolve{}’s role as a foundational component in scientific modeling across disciplines.

\section{Mathematical Description}\label{sec:mathematical_description}

This section introduces the mathematical framework for numerically solving nonlinear problems and demonstrates the built-in support for such problems in \nlsolve. A nonlinear problem is defined as:
\begin{equation}
    \text{Find } u^\ast \text{ s.t. } ~~~ \func{f}{u^\ast, \theta} = 0 \label{eq:nonlinear-root-finding}
\end{equation}
where $f : \mathbb{R}^n \times \mathbb{P} \to \mathbb{R}^n$ and $\theta \in \mathbb{P}$ is a set of fixed parameters. In this paper, we will assume that $f$ is continuous and differentiable\footnote{For Halley's method, we additionally assume twice-differentiability.}. For a more rigorous and detailed discussion on numerical solvers for nonlinear systems, we refer the readers to \cite[Chapter 11]{wright1999numerical}.

\subsection{Numerical Algorithms for Nonlinear Equations}\label{subsec:numerical-root-finding}

Numerical nonlinear solvers, alternatively referred to as iterative algorithms, start with an initial guess $u_0$ and iteratively refine it until a convergence criterion is satisfied, typically, $\norm{\func{f}{u_k, \theta}}_\infty \leq \texttt{abstol}$. (where $\norm{\cdot}_\infty$ denotes the max-norm). Newton-Raphson is a powerful iterative technique for finding the roots of a continuously differentiable function $f(u, \theta)$. It iteratively refines an initial guess $u_0$ using the formula:
\begin{equation}
    \begin{aligned}
        \mathcal{J}_k \delta u_\text{newton} & = -f(u_k, \theta) \label{eq:newton_descent} \\
        u_{k + 1}                            & = u_k + \delta u_\text{newton},
    \end{aligned}
\end{equation}
where $\mathcal{J}_k$ is the Jacobian of $f(u, \theta)$ with respect to $u$, evaluated at $u_k$. This method exhibits rapid convergence~\cite[Theorem~11.2]{wright1999numerical} when the initial guess is sufficiently close to a root. Furthermore, it requires only the function and its Jacobian, making it computationally efficient for many practical applications. Halley's method enhances the Newton-Raphson method, leveraging information from the second total derivative of the function to achieve cubic convergence. It refines the initial guess $u_0$ using:
\begin{equation}
    \begin{aligned}
        a^\nu_k                & = \delta u_\text{newton}   &
        \mathcal{J}_k b^\nu_k  & = \mathcal{H}_k a^\nu_k  a^\nu_k                                               \\
        \delta u_\text{halley} & = \left( a^\nu_k \circ a^\nu_k \right) \oslash \left( a^\nu_k + \frac{1}{2} b^\nu_k \right) &
        u_{k + 1}              & = u_k + \delta u_\text{halley}
    \end{aligned}
\end{equation}
where $\mathcal{H}_k = \nabla_u^2 f(u, \theta)$ is the second total derivative of $f(u, \theta)$ w.r.t. $u$, evaluated at $u_k$, and $\circ$ and $\oslash$ denote element-wise multiplication and division respectively. Despite the cubic convergence, the additional cost of computing $\mathcal{H}_k$ with naive forward or reverse mode AD~\cite{bettencourt2019taylor} often makes Halley's method practically infeasible. \citet{tan2025scalable} uses Taylor-mode AD to implement higher-order solvers in \nlsolve, highlighting the flexibility of its modular design.

Another widely used solver for nonlinear problems is the Levenberg-Marquardt algorithm, which interpolates between gradient descent and the Gauss-Newton method by introducing a damping term~\cite{transtrum2012improvements}. It is especially effective when the Jacobian is ill-conditioned or when the initial guess is far from a root. To further improve convergence in narrow, curved valleys of the cost surface, \citet{transtrum2012improvements} introduces a geodesic acceleration correction -- a second-order update along the current search direction that accounts for local curvature. This correction is inexpensive to compute and improves robustness without requiring full second derivatives. 

Multi-step methods are a special class of nonlinear solvers that take multiple intermediate steps between $u_k$ and $u_{k + 1}$ to obtain a higher order of convergence. For example, Potra \& Pt\'{a}k~\cite{potra1980nondiscrete} uses the same Jacobian twice to compute the final descent direction.
\begin{equation}
    \begin{aligned}
        \mathcal{J}_k \delta u_1 & = -f(u_k, \theta)       & y_k       & = u_k + \delta u_1 \\
        \mathcal{J}_k \delta u_2 & = -f(y_k, \theta)       & u_{k + 1} & = y_k + \delta u_2.
    \end{aligned}
\end{equation}
This method provides higher convergence orders without reliance on higher-order derivatives. \cite{singh2023simple} summarizes other multi-step schemes that provide higher-order convergence using only first-order derivatives. All the methods we discussed in this section are local algorithms, and their convergence relies on having a good initial guess. We will discuss some techniques to facilitate the global convergence of these methods in the following section.

\subsection{Globalization Strategies}\label{subsec:globalization-strategies}

Globalization strategies allow us to enhance the local algorithms described in \Cref{subsec:numerical-root-finding} and ensure global convergence under suitable assumptions~\cite[Theorem 3.2, 4.5, 4.6]{wright1999numerical}.

\subsubsection{Line Search}

\begin{figure}[t]
    \begin{minipage}[c]{0.4\textwidth}
        \centering
        \begin{minted}[breaklines,escapeinside=||,mathescape=true, linenos, numbersep=3pt, gobble=2, frame=lines, fontsize=\ssmall, framesep=2mm]{julia}
using NonlinearSolve, LineSearches

function generalized_rosenbrock(x, p)
    return vcat(1 - x[1], 10 .* (x[2:end] .- x[1:(end - 1)] .* x[1:(end - 1)]))
end

x_start = vcat(-1.2, ones(9))

prob = NonlinearProblem(generalized_rosenbrock, x_start; abstol = 1e-8)

solve(prob, NewtonRaphson(; linesearch = <LINE SEARCH>))
        \end{minted}
    \end{minipage}
    \hfill
    \begin{minipage}[c]{.59\textwidth}
        \centering
        \includegraphics[width=\textwidth]{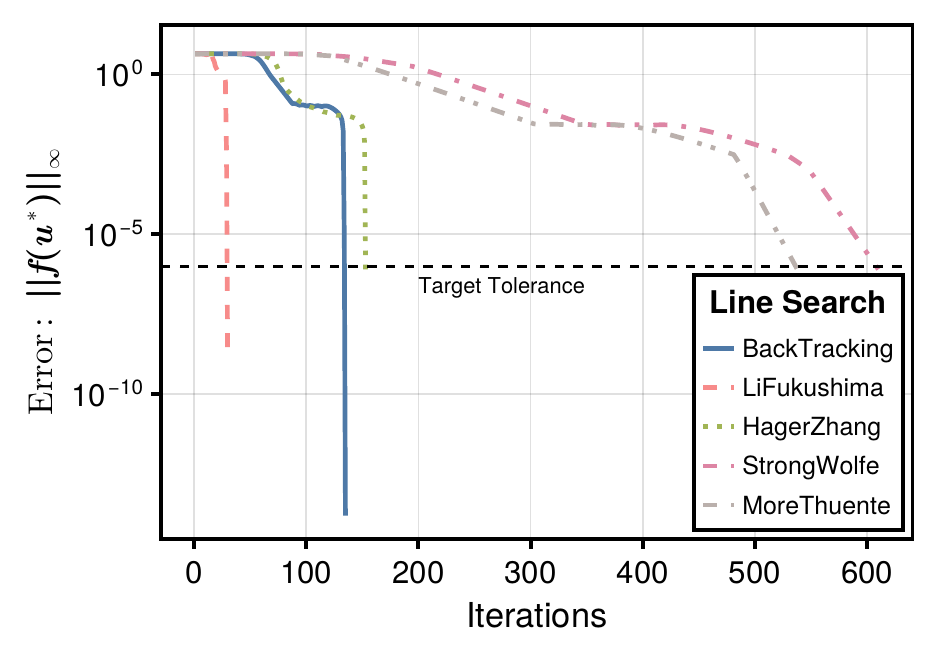}
    \end{minipage}
    \vspace{-1em}
    \caption{\nlsolve{} allows seamless switching between different line search routines. We can specify \texttt{HagerZhang()}, \texttt{BackTracking()}, \texttt{MoreTheunte()}, \texttt{StrongWolfe()}, or \texttt{LiFukushimaLineSearch()} to opt-in to using line search. Line search enables \texttt{NewtonRaphson} to converge faster on the generalized Rosenbrock problem [\Cref{eq:gen_rosenbrock_function}], while it fails to converge without a line search.}\label{fig:nr-line-search}
    \vspace{-2em}
\end{figure}

Consider the $N$-dimensional generalized Rosenbrock function, which is defined as:

\begin{align}
    f_i(x, \theta) = \begin{cases}
        1 - x_1 & \text{if } i = 1\\
        10  \left(x_{i} - x_{i - 1}^2\right) & \text{if } i \in \lbrace 2, \dots, N\rbrace
    \end{cases}\label{eq:gen_rosenbrock_function}
\end{align}
If we initialize the problem with $(u_0)_1 = -1.2$, Newton-Raphson fails to converge to the solution [\Cref{fig:nr-line-search}]. Line search methods are globalization strategies that avoid such failures by determining a positive step length $\alpha_k$ given the current iterate $u_k$ and the search direction $\delta u_k$
\begin{equation}
    u_{k + 1} = u_k + \alpha_k \delta u_k.
\end{equation}
The direction is often given by the Newton descent direction, Steepest descent, or one of the Multi-Step Schemes described in \Cref{subsec:numerical-root-finding}. The optimal choice for the step length for a given merit function $\phi(u_k, \alpha) = \frac{1}{2} \norm{f(u_k + \alpha_k \delta u_k, \theta)}^2$ is given by:
\begin{equation}
    \alpha^* \in \argmin_{\alpha > 0}  \phi(u_k, \alpha)\label{eq:exact-line-search}
\end{equation}
However, solving a global optimization problem on each step of the iterative nonlinear solver is prohibitively expensive. Instead, line search methods rely on selecting a set of candidate step sizes and terminating the search based on certain conditions:

\paragraph{Armijo Rule and Curvature Criteria} The Armijo Rule or Sufficient Decrease criteria states that $\alpha_k$ is acceptable only if there is a sufficient decrease in the objective function:
\begin{align}
    \phi(u_k, \alpha) &\leq \phi(u_k, 0) + c_1 \alpha \nabla_{u_k} \phi(u_k, 0)^T \delta u_k \\
    &= \phi(u_k, 0) + c_1 \alpha \underbrace{  f(u_k, \theta)^\top \mathcal{J}_k \delta u_k }_{\substack{\text{directional derivative} \\ \text{of } \frac{1}{2}(f)^2 \text{ at } u_k \text{ along } \delta u_k}}  \label{eq:sufficient-decrease}
\end{align}
Additionally, we use the curvature condition to filter out values for $\alpha_k$ that satisfy sufficient decrease but are very close to the current iterate. This ensures that the algorithm makes reasonable progress at every step.
\begin{equation}
    \nabla_{u_k} \phi(u_k, \alpha_k)^T \delta u_k \geq c_2 \nabla_{u_k} \phi(u_k, 0)^T \delta u_k
\end{equation}
where $c_2 \in (c_1, 1)$ and $c_1 \in (0, 1)$. These two conditions are collectively known as the Wolfe Conditions~\cite{wolfe1969convergence, wolfe1971convergence}.

\paragraph{Strong Wolfe Condition} More sophisticated line search methods also enforce the curvature condition and its strong variant~\cite{more1994line}, which bounds the absolute value of the directional derivative to avoid steps that are too aggressive in either direction
\begin{align}
     & \phi(u_k, \alpha) \leq \phi(u_k, 0) + c_1 \alpha \nabla_{u_k} \phi(u_k, 0)^T \delta u_k                                   \\
     & \lvert \nabla_{u_k} \phi(u_k, \alpha_k)^T \delta u_k \rvert \geq c_2 \lvert \nabla_{u_k} \phi(u_k, 0)^T \delta u_k \rvert
\end{align}

In practice, line search methods are implemented with approximations and heuristics to balance robustness and computational efficiency. \texttt{Backtracking Line Search} is one such simple strategy, which starts with an initial step length (often $\alpha_0 = 1$) and is iteratively reduced by a factor $\rho \in (0, 1)$ until the sufficient decrease condition [\Cref{eq:sufficient-decrease}] is satisfied. This ensures that the algorithm makes consistent progress, assuming the search direction $\delta u_k$ is a descent direction for the merit function $\phi(u_k, \alpha)$. Because sufficient decrease is explicitly enforced, backtracking line search typically does not require checking the curvature condition. A detailed description can be found in~\cite[Section 3.1]{wright1999numerical}.
The \texttt{LineSearches.jl} package~\cite{Optim.jl-2018} provides a flexible suite of line search algorithms, including \texttt{BackTracking}~\cite{wright1999numerical}, \texttt{HagerZhang}~\cite{hager2006algorithm}, and \texttt{MoreThuente}~\cite{more1994line}, each with different strategies for satisfying the Wolfe or Strong Wolfe conditions. These are supported in \nlsolve{} through a uniform API that enables easy experimentation. In \Cref{fig:nr-line-search}, we demonstrate how line search methods can resolve convergence issues in Newton-Raphson when initialized far from the solution, using the generalized Rosenbrock function [\Cref{eq:gen_rosenbrock_function}].

\begin{figure}
    \centering
    \includeinkscape[width=\textwidth]{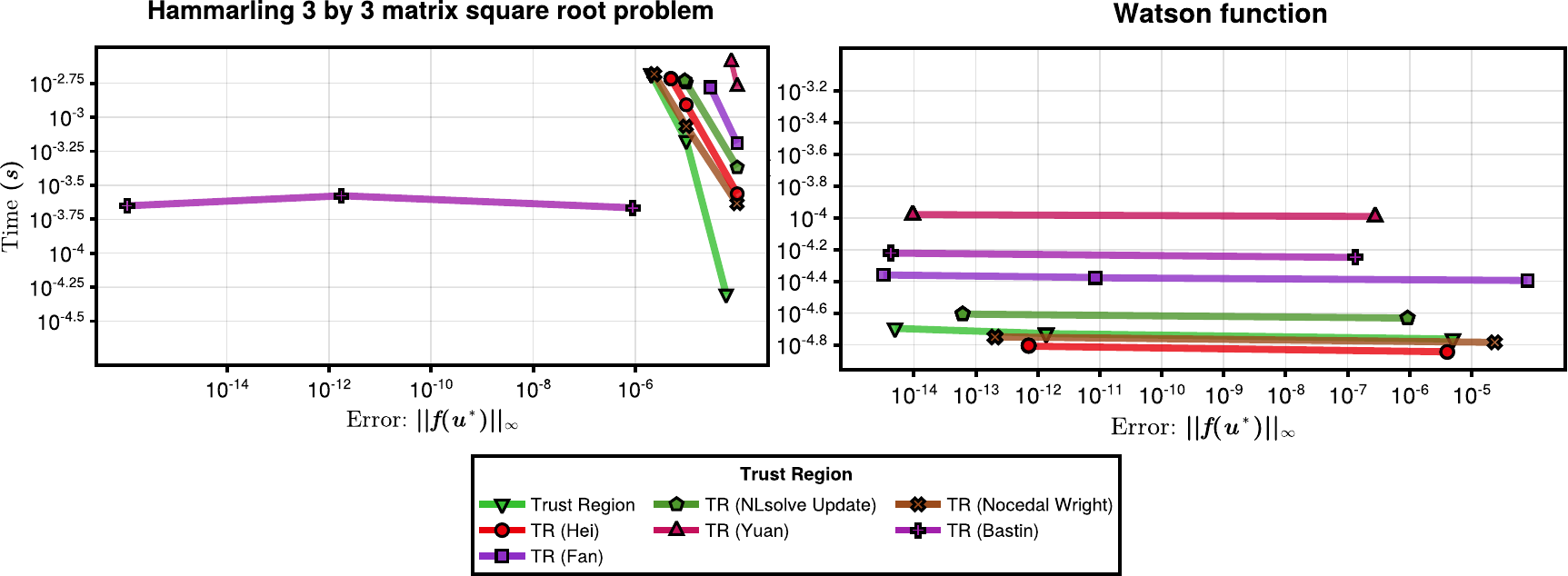}
    \caption{\textbf{Work-Precision Diagrams for Two Instances of the 23 Test Problems}: Trust-region radius update schemes directly affect the convergence of the algorithm, \nlsolve{} implements various update schemes and allows switching them with a unified API.}\label{fig:trust-region-rus-effect}
    \vspace{-2em}
\end{figure}

\subsubsection{Trust-Region Methods}

As an alternative to using line search as a globalization strategy, trust-region methods restrict the next iterate to lie in a local region around the current iterate. These methods estimate a quadratic model for the objective function:
\begin{equation}
    \func{m_k}{\delta u_k} = \frac{1}{2}\norm{f(x, \theta)}^2 + \delta u_k^T \mathcal{J}_k^T f(x, \theta) + \frac{1}{2} \delta u_k^T \mathcal{J}^T_k \mathcal{J}_k \delta u_k
\end{equation}
and attempt to solve the following sub-problem:
\begin{equation}
    \delta u_k^\ast \in ~\argmin_{\delta u_k}~ \func{m_k}{\delta u_k} \qquad \text{s.t. } \norm{\delta u_k} \leq \Delta_k\label{eq:tr_subproblem}
\end{equation}
where $\Delta_k$ is the trust-region radius at the current iterate $k$. To perform an update to the current iterate and the trust-region radius, we compute the ratio of the actual reduction to the predicted reduction ($\rho_k$):
\begin{equation}
    \rho_k = \frac{\norm{f(x, \theta)}^2 - \norm{f(x + \delta u_k, \theta) }^2}{\norm{f(x, \theta)}^2 - \norm{f(x, \theta) + \mathcal{J}_k \delta u_k}^2}
\end{equation}
If $\rho_k$ is greater than a threshold, i.e., there is a strong agreement between the predicted model and true system, we accept the step and increase the trust-region radius $\Delta_k$. If $0 < \rho_k << 1$, we accept the step but keep the trust-region radius unaffected. Otherwise, the trust-region is shrunk, and the step is rejected. The exact specification for thresholds for these choices depends on the underlying algorithms being used. \nlsolve{} implements additional sophisticated algorithms to update the trust-region -- \texttt{Hei}~\cite{hei2003self}, \texttt{Yuan}~\cite{yuan2015recent}, \texttt{Fan}~\cite{fan2006convergence}, \texttt{NocedalWright}~\cite{wright1999numerical}, \texttt{NLsolve}~\cite{Optim.jl-2018}, and \texttt{Bastin}~\cite{bastin2010retrospective}. In \Cref{fig:trust-region-rus-effect} and \Cref{sec:results}, we demonstrate how different radius update schemes achieve a compromise between the computational cost and solution accuracy.

The trust-region subproblem [\Cref{eq:tr_subproblem}] can be solved approximately using Powell's Dogleg Method \cite{powell1970hybrid}, which combines the descent directions from Newton Descent ($\delta u_\text{newton} = -\mathcal{J}^{-1} f(u, \theta)$) and the Steepest Descent ($\delta u_\text{sd} = -\mathcal{J}^T f(u, \theta)$). The algorithm proceeds by computing $\delta u_\text{newton}$, and takes the full step if $\norm{\delta u_\text{newton}} \leq \Delta$. Otherwise, it computes the unconstrained minimizer of the model along the steepest descent direction. If the minimizer is outside the trust-region, it accepts the descent vector's intersection with the hypersphere of radius $\Delta$ if $\norm{\delta u_\text{sd}} \geq \Delta$. Finally, if none of these conditions hold, we compute the interpolated vector $\delta u = (1 - \tau) \delta u_\text{sd} + \tau \delta u_\text{newton}$, where $\tau$ is the largest value $ \in (0, 1)$ s.t. $\norm{\delta u} \leq \Delta$.

\subsection{Sensitivity Analysis}\label{subsec:sensitivity_analysis}

Sensitivity analysis for nonlinear systems is an integral part of several applications like differentiable physics simulators, deep equilibrium models~\cite{baideep2019, pal2023continuous, pal2023efficient, pal2022lux}, and optimizing design parameters~\cite{ruppert2023adjoint}. Computing gradients for these systems using naive automatic differentiation requires differentiating the nonlinear solve process, which is computationally expensive since it increases the order of differentiation for the solver. Instead, we implement a black-box adjoint method using the implicit function theorem that differentiates the solution directly, independent of the exact solver being used. For a continuous differentiable function $f: \mathbb{R}^n \times \mathbb{P} \rightarrow \mathbb{R}^n$ and $u^\ast \in \mathbb{R}^n$ and $\theta \in \mathbb{P}$ s.t.:
\begin{equation}
    \func{f}{u^\ast, \theta} = 0
\end{equation}
and the Jacobian w.r.t $u^\ast$ being non-singular, the Implicit Function Theorem (IFT) provides an efficient way to compute the Jacobian~\cite{blondel2021opez, johnson2012notes}:
\begin{equation}
     \frac{\mathrm{d}u^\ast}{\mathrm{d}\theta} = \left( \nabla_u f(u^\ast, \theta) \right)^{-1} \nabla_\theta f(u^\ast, \theta)
\end{equation}
\nlsolve{} has specialized Krylov methods to compute the sensitivity efficiently without constructing these Jacobians explicitly using vector-Jacobian Products (VJP). Consider the following example, where we use the forward mode AD framework \texttt{ForwardDiff.jl}~\cite{RevelsLubinPapamarkou2016} and reverse mode AD framework \texttt{Zygote.jl}~\cite{Zygote.jl-2018, innes2018fashionable} to differentiate through the solution of $f(u, \theta) = u^2 - \theta = 0$ with $g(u^\ast, \theta) = \sum_i (u_i^\ast)^2$:

\begin{minted}[breaklines,escapeinside=||,mathescape=true, linenos, numbersep=3pt, gobble=2, frame=lines, fontsize=\ssmall, framesep=2mm]{julia}
using NonlinearSolve, SciMLSensitivity, ForwardDiff, Zygote

f(x, p) = @. x^2 - p
g(u_star) = sum(abs2, u_star)

solve_nlprob(p) = g(solve(NonlinearProblem(f, [1.0, 2.0], p)).u)

p = [2.0, 5.0]

ForwardDiff.gradient(solve_nlprob, p)
Zygote.gradient(solve_nlprob, p)
\end{minted}

\subsection{Matrix Coloring \& Sparse Automatic Differentiation}\label{subsec:sparse_colored_ad}

The Jacobian computation and linear solver are typical bottlenecks for numerical nonlinear equation solvers. If the sparsity pattern of a Jacobian is known, then computing the Jacobian can be done with much higher efficiency~\cite{averick1994computing}. Sparsity patterns for arbitrary programs can be automatically generated using numerical techniques~\cite{giering2005generating, walther2009getting} or symbolic methods~\cite{gowda2019sparsity}. Given the sparsity pattern, we can use graph coloring algorithms~\cite{gebremedhin2005color, sparsediff} to compute the matrix colors for the given sparse matrix.

In \nlsolve, we can compute the sparse Jacobian using both forward and reverse mode AD. AD tools compute the directional derivative using Jacobian-vector products (JVPs) $\mathcal{J}  u$ for forward mode AD and VJPs $\mathcal{J}^T  v$ for reverse mode AD. To compute a dense $m \times n$ Jacobian, forward mode AD (and finite differencing) computes $n$ JVPs, and reverse mode AD computes $m$ VJPs using the standard basis vectors. Forward mode AD fills in the Jacobian one column at a time, and reverse mode AD fills in the Jacobian one row at a time. We can chunk multiple JVPs or VJPs together for Jacobians with known sparsity patterns to reduce the overall computation cost. Consider the following sparsity pattern for a Jacobian $\mathcal{J}_{\text{sparse}}$ where $\bullet$ denotes the non-zero elements of the Jacobian. To use forward mode AD, we need to compute the column coloring, and for reverse mode we need the row coloring (which can be obtained by coloring the transposed Jacobian).
\begin{equation}
    \mathcal{J}_{\text{sparse}} = \begin{bmatrix}
        \bullet &         &         &         &         \\
                & \bullet & \bullet &         &         \\
                &         &         & \bullet &         \\
        \bullet & \bullet &         &         & \bullet \\
                &         &         &         & \bullet
    \end{bmatrix} \qquad  \mathcal{J}^{(\text{col})}_{\text{sparse}} = \begin{bmatrix}
        \color{red}{\blacktriangleright} &                            &                                  &                                  &                              \\
                                         & \color{blue}{\blacksquare} & \color{red}{\blacktriangleright} &                                  &                              \\
                                         &                            &                                  & \color{red}{\blacktriangleright} &                              \\
        \color{red}{\blacktriangleright} & \color{blue}{\blacksquare} &                                  &                                  & \color{green}{\blacklozenge} \\
                                         &                            &                                  &                                  & \color{green}{\blacklozenge}
    \end{bmatrix} \qquad \mathcal{J}^{(\text{row})}_{\text{sparse}} = \begin{bmatrix}
        \color{blue}{\blacksquare}   &                              &                            &                            &                              \\
                                     & \color{blue}{\blacksquare}   & \color{blue}{\blacksquare} &                            &                              \\
                                     &                              &                            & \color{blue}{\blacksquare} &                              \\
        \color{green}{\blacklozenge} & \color{green}{\blacklozenge} &                            &                            & \color{green}{\blacklozenge} \\
                                     &                              &                            &                            & \color{blue}{\blacksquare}
    \end{bmatrix}
\end{equation}
Matrix Coloring~\cite{gebremedhin2005color} allows us to avoid perturbation confusion when we perturb multiple rows/columns of the same color together. For example, in forward mode AD perturbing the columns with color $\color{red}{\blacktriangleright}$, we get $\left[u_1 + \varepsilon, u_2, u_3 + \varepsilon, u_4 + \varepsilon, u_5\right]$ ($a + \varepsilon$ denotes a Dual number). Propagating this vector, we obtain the result $\left[v_1 + \delta_1\varepsilon, \dots, v_5 + \delta_5\varepsilon\right]$. Since all the columns of the Jacobian with the same color are non-overlapping, we can uniquely map the partials to the Jacobian entries: $\left(\mathcal{J}_{\text{sparse}}\right)_{11} = \delta_1$, $\left(\mathcal{J}_{\text{sparse}}\right)_{23} = \delta_2$,
$\left(\mathcal{J}_{\text{sparse}}\right)_{34} = \delta_3$,
$\left(\mathcal{J}_{\text{sparse}}\right)_{41} = \delta_4$, and from the sparsity pattern $\delta_5 = 0$. For reverse mode, we proceed similarly, except the perturbation is now on the output. Hence, to compute $\mathcal{J}^{(\text{col})}_{\text{sparse}}$, we need to perform 3 JVPs (once each for $\color{red}{\blacktriangleright}$, $\color{blue}{\blacksquare}$, and $\color{green}{\blacklozenge}$) compared to 5 JVPs for a $5 \times 5$ dense Jacobian. Similarly, since reverse mode materializes the Jacobian one row at a time, we need 2 VJPs (once each for $\color{blue}{\blacksquare}$, and $\color{green}{\blacklozenge}$) compared to 5 VJPs for the dense counterpart. With this simple example, we demonstrate a specific case where sparse reverse mode is more efficient than sparse forward mode and vice versa. \nlsolve{} provides abstractions, allowing seamless switching between the modes for maximum performance.  In \Cref{subsec:automatic-sparsity-detection}, we benchmark the scaling of sparse AD with regular AD on a steady-state PDE.

\begin{figure}[t]
    \centering
    \includegraphics[width=\textwidth]{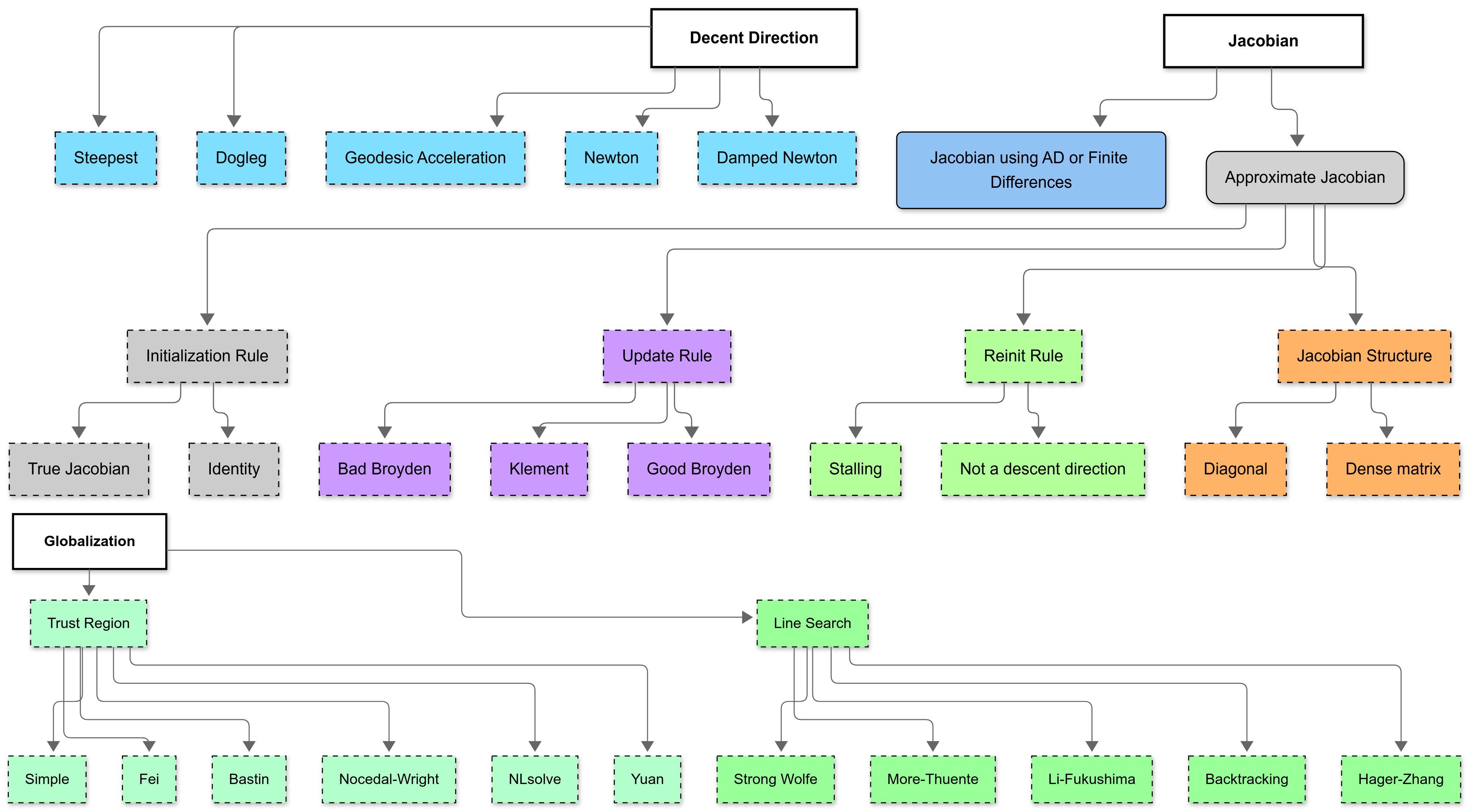}
    \caption{\nlsolve{} has a modular architecture with the core building blocks: Jacobian Computation Algorithm, Globalization Strategy, and Descent Direction. By combining these fundamental components, \nlsolve{} facilitates creating diverse and powerful solver algorithms tailored to specific problem characteristics. This composability feature underscores the framework's versatility in tackling a wide range of nonlinear equations, showcasing the potential for innovative algorithm development through the flexible integration of different strategies and techniques.}\label{fig:solver_options}
    \vspace{-2em}
\end{figure}

\section{Special Capabilities}\label{sec:special-capabilities}

In the previous section, we provided an overview of solving nonlinear equations. Here, we will focus on specific capabilities of our package \nlsolve. Let's cover the typical workflow of using \nlsolve. First, we define the nonlinear equations. We will use a mutating function to compute our residual value from $u$ and $p$ and store it in $F$.
\begin{minted}[breaklines,escapeinside=||,mathescape=true, numbersep=3pt, gobble=2, frame=lines, fontsize=\ssmall, framesep=2mm]{julia}
function f!(F, x, p)
    a, b, c = p
    F[1] = (x[1] + a) * (x[2]^3 - b) + c
    F[2] = sin(x[2] * exp(x[1]) - 1)
    return nothing
end
\end{minted}
In Julia, mutating function names are appended with a \texttt{!} as a convention. Next, we construct a \texttt{NonlinearFunction} -- a wrapper over \texttt{f!} that can store additional information like user-defined Jacobian functions and sparsity patterns.
\begin{minted}[breaklines,escapeinside=||,mathescape=true, numbersep=3pt, gobble=2, frame=lines, fontsize=\ssmall, framesep=2mm]{julia}
nlfunc = NonlinearFunction(f!)
\end{minted}
Next, we construct the \texttt{NonlinearProblem} with 3 arguments -- the nonlinear function \texttt{nlfunc}, the initial guess $u_0 = [0.0, 0.0]$, and the fixed parameters $\theta = [3.0, 7.0, 18.0]$.
\begin{minted}[breaklines,escapeinside=||,mathescape=true, numbersep=3pt, gobble=2, frame=lines, fontsize=\ssmall, framesep=2mm]{julia}
prob = NonlinearProblem(nlfunc, [0.0, 0.0], [3.0, 7.0, 18.0])
\end{minted}
Finally, calling \texttt{solve} on the problem uses a poly-algorithm [\Cref{subsec:poly-algorithm-defaults}] to solve this problem. To use any other algorithm, we can pass it as the second argument.
\begin{minted}[breaklines,escapeinside=||,mathescape=true, numbersep=3pt, gobble=2, frame=lines, fontsize=\ssmall, framesep=2mm]{julia}
solve(prob)
solve(prob, NewtonRaphson())
\end{minted}

\subsection{Composable Building Blocks}\label{subsec:composable-blocks}

\nlsolve{} stands out for its modular design, allowing users to combine various building blocks to create custom solver algorithms [\Cref{fig:solver_options}]. This section delves into the flexibility and power of \nlsolve's architecture, focusing on its core components: Jacobian Computation Algorithm, Globalization Strategy, and Descent Direction.

\begin{table}[tbp]
   \centering
   \begin{adjustbox}{width=\textwidth,center}
       \begin{tabular}{l l l l l l l l l l}
            \toprule
            \multirow{3}{*}{\textbf{Algorithm}} & & \multicolumn{5}{c}{\textbf{Jacobian Strategy}} &  & \multirow{3}{*}{\textbf{Globalization}} & \multirow{3}{*}{\textbf{Descent}}\\
            \cmidrule{3-7}
            & &  \textbf{AD / FD} & \multirow{2}{*}{\textbf{Structure}} & \textbf{Init.} & \textbf{Reinit} & \textbf{Update}  & & & \\
            &  & \textbf{Jac.} & & \textbf{Rule} & \textbf{Rule} & \textbf{Rule} &  & & \\
            \midrule
            \texttt{NewtonRaphson} & & \cmark & - & - & - & - & & - & \texttt{Newton}\\
            \texttt{NR (BackTracking)} & & \cmark & - & - & - & - & & \texttt{BackTracking} & \texttt{Newton}\\
            \texttt{Trust-Region} & & \cmark & - & - & - & - & & \texttt{Simple} & \texttt{Dogleg}\\
            \texttt{TR (Bastin)} & & \cmark & - & - & - & - & & \texttt{Bastin} & \texttt{Dogleg}\\
            & & & & & & & & & \\
            \texttt{Broyden} & & \xmark & \texttt{Dense} & \texttt{Identity} & & & & & \\
            \texttt{Limited Mem.} & & \multirow{2}{*}{\xmark} & \texttt{Low Rank} & \multirow{2}{*}{\texttt{Identity}} & \texttt{Not a Des.} & \multirow{2}{*}{\texttt{Good Broyden}} &  & \multirow{2}{*}{-} & \multirow{2}{*}{\texttt{Newton}}\\
            \texttt{Broyden} & & & \texttt{Matrix} & & \texttt{Direction} & & & & \\
            \texttt{Mod. Broyden} & & \xmark & \texttt{Dense} & \texttt{True Jac.} & &   & & &\\
            & & & & & & & & & \\
            \texttt{Klement} & & \xmark & \texttt{Diagonal} & \texttt{Identity} & \texttt{Stalling} & \texttt{Klement} & & - & \texttt{Newton}\\
            \texttt{PseudoTransient} & & \cmark & - & - & - & - & & - & \texttt{DampedNewton}\\
            \multirow{2}{*}{\texttt{LevenbergMarquart}} & & \multirow{2}{*}{\cmark} & \multirow{2}{*}{-} & \multirow{2}{*}{-} & \multirow{2}{*}{-} & \multirow{2}{*}{-} & & \texttt{Indirect TR} & \multirow{2}{*}{\texttt{DampedNewton}}\\
            & & & & & & & & \texttt{Update} &  \\
            \bottomrule
        \end{tabular}
    \end{adjustbox}
    \vspace{1em}
    \caption{\textbf{Building Blocks Provided by \nlsolve{} to Construct Nonlinear Root-Finding Algorithms}}\label{tab:building-blocks}
    \vspace{-2em}
\end{table}
\begin{enumerate}
    \item \ul{Jacobian Computation}: \nlsolve{} supports a range of methods for Jacobian computation. Firstly, we support various automatic differentiation and finite differencing backends to compute the Jacobian. Additionally, to support Quasi-Newton methods, we provide different abstractions to iteratively approximate the Jacobian.

    \item \ul{Globalization Strategy}: This module facilitates the global convergence of the local solvers [\Cref{subsec:globalization-strategies}] and can seamlessly integrate with the other blocks.

    \item \ul{Descent Direction}: These determine the direction vector given the Jacobian (or higher-order derivatives) and the current value.
\end{enumerate}

\Cref{tab:building-blocks} demonstrates the building blocks used in \nlsolve{} to construct numerical nonlinear root-finding problems. The modular architecture of \nlsolve{} supports a broad range of algorithms and empowers users to craft sophisticated solver combinations. Among the possibilities, users can integrate Pseudo-Transient methods with globalization strategies or pair Quasi-Newton methods with globalization techniques. This flexibility opens the door to creating tailored solutions that can more effectively navigate challenging problem landscapes, highlighting \nlsolve's capability to adapt and optimize across a diverse set of nonlinear equation challenges.

\subsection{Smart Poly-Algortihm Defaults}\label{subsec:poly-algorithm-defaults}

\begin{figure}[t]
    \centering
    \includegraphics[width=\textwidth]{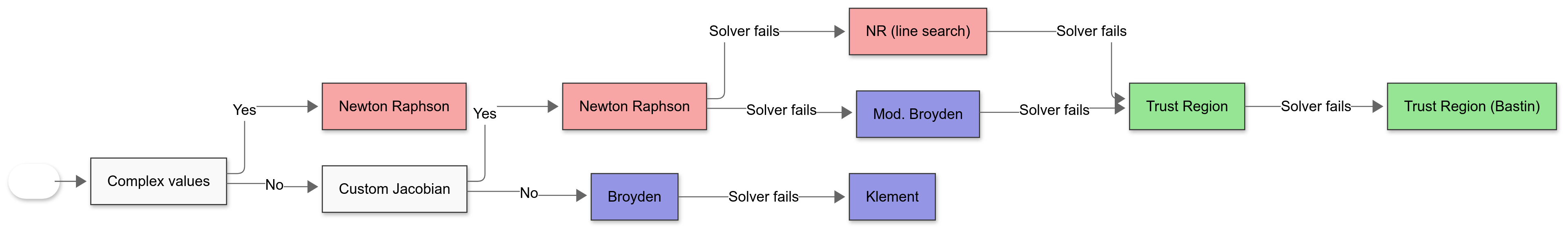}
    \caption{\textbf{Default Nonlinear Solve Poly-Algorithm} focused on balancing speed and robustness. \nlsolve{} first tries less robust Quasi-Newton methods for better performance and then tries more robust techniques if the faster ones fail. Based on the choice of linear solver, these methods automatically switch to a Jacobian-free variant.}\label{fig:nonlinear_solve_default_alg}
\end{figure}

\begin{figure}[t]
    \centering
    \includegraphics[width=\textwidth]{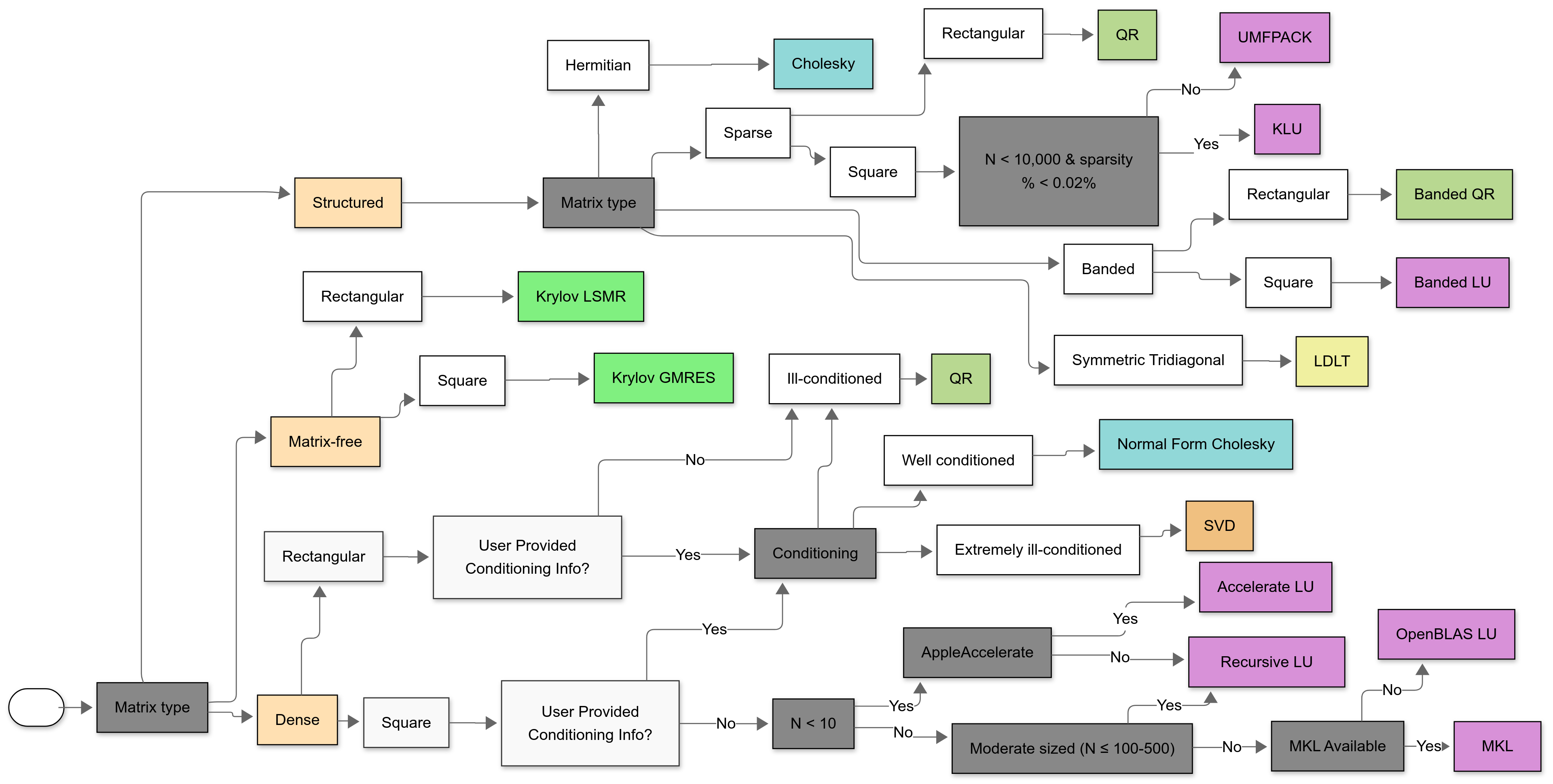}
    \caption{\textbf{Decision Flowchart for Default Linear Solver Selection}: Starting with the determination of the Jacobian type (e.g., square, rectangular, symmetric, structured), the flowchart guides through various conditions such as sparsity and conditioning to select the most appropriate solver. Options range from specialized solvers like KLU for sparse matrices to robust methods like SVD for extremely ill-conditioned cases. This adaptive approach ensures that \nlsolve{} utilizes the most efficient linear solver, enhancing performance and accuracy for a wide array of nonlinear problems.}\label{fig:linear_solve_default_alg}
    \vspace{-2em}
\end{figure}

\nlsolve{} comes with a set of well-benchmarked defaults designed to be fast and robust in the average case. For nonlinear problems with roots in pre-specified intervals $\left[a, b\right]$, we default to Interpolate, Truncate, and Project (\texttt{ITP})~\cite{oliveira2020enhancement}. For more advanced problems, we default to a poly-algorithm that selects the internal solvers based on the problem specification. Our selection happens two-fold, first at the nonlinear solver level [\Cref{fig:nonlinear_solve_default_alg}] and next at the linear solver level\footnote{Default linear solver selection happens in \texttt{LinearSolve.jl}.} [\Cref{fig:linear_solve_default_alg}].

\paragraph{Default Nonlinear Solver Selection}

By default, we use a poly-algorithm designed to balance speed and robustness. The solver begins with Quasi-Newton methods -- \texttt{Broyden}~\cite{broyden1965class} and \texttt{Klement}~\cite{klement2014using}. If these fail, it falls back to a modified Broyden method that initializes the Jacobian with $\nabla_u f(u, \theta)$. If all Quasi-Newton methods fail, the algorithm switches to Newton-Raphson, then to Newton-Raphson with line search, and finally to trust-region methods. When a custom Jacobian is provided or the system has 25 or fewer states, we skip the Quasi-Newton phase and directly apply first-order methods. This layered strategy ensures robustness by default while allowing users to opt into faster alternatives. %

\paragraph{Default Linear Solver Selection}

We follow a structured decision flow (see \Cref{fig:linear_solve_default_alg}) to select the most appropriate linear solver for systems of the form $Ax = b$. The decision process begins with determining the structure and shape of the Jacobian -- e.g., whether it's structured, symmetric, square, or rectangular. For structured or symmetric matrices, we prefer specialized solvers like \texttt{Cholesky}, \texttt{LDLT}, or \texttt{NormalFormCholesky}. For general square matrices, we apply \texttt{LU} or \texttt{QR}, with variants chosen based on matrix size and sparsity. If the system is small and dense, we default to native Julia LU routines~\cite{toledo1997locality, recursivefactorization}. For larger matrices, we use OpenBLAS or MKL-backed LU routines if available.

For sparse matrices with $N > 10,000$ and density $< 0.01$, we use solvers like \texttt{KLU} or \texttt{UMFPACK}. When users provide a custom preconditioner, we use iterative solvers such as \texttt{GMRES}~\cite{saad1986gmres, montoison-orban-2023}. If the matrix is matrix-free, we automatically invoke iterative Krylov methods without materializing the Jacobian. For ill-conditioned or nearly singular systems, we fall back to robust solvers like \texttt{QR}. In extremely ill-conditioned cases, we default to \texttt{SVD} to ensure numerical stability. This adaptive selection strategy ensures that \nlsolve{} consistently selects the most efficient and stable solver for a wide range of nonlinear problems.

\begin{figure}[t]
    \centering
    \adjustbox{trim=4em 3em 3em 3em, clip}{
        \includeinkscape[width=1.15\textwidth]{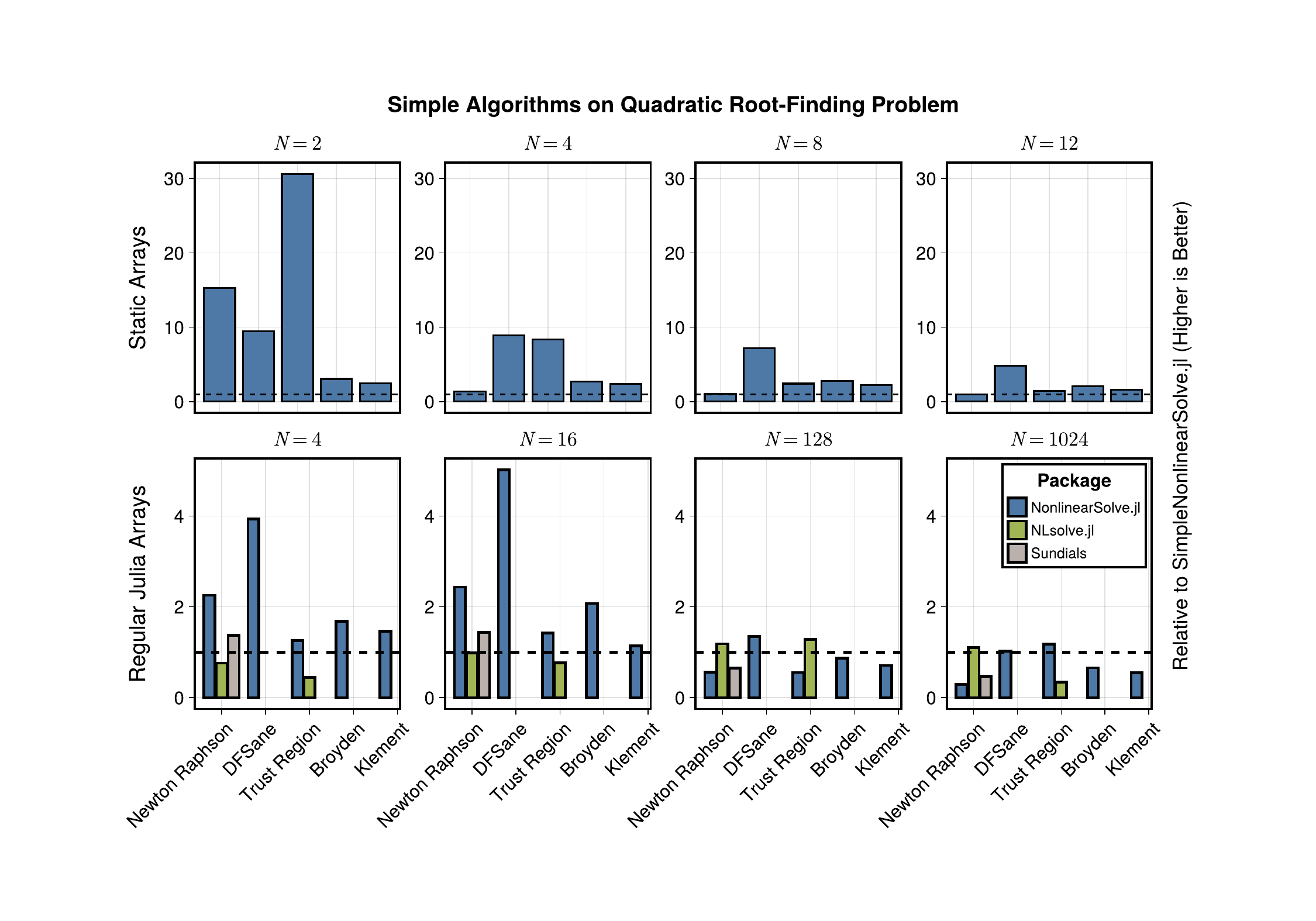}
    }
    \caption{\texttt{SimpleNonlinearSolve.jl} provides the fastest algorithms for StaticArray Problems, outperforming the same algorithms from \nlsolve{} by a significant margin. For regular arrays, \texttt{SimpleNonlinearSolve.jl} is competitive to standard solvers up to a certain problem size threshold.}\label{fig:quadratic_problem_simplenonlinearsolve}
\end{figure}

\subsection{Non-Allocating Static Algorithms inside GPU Kernels}\label{subsec:simplenonlinearsolve}

Julia has long supported general-purpose GPU programming (GPGPU) through \texttt{CUDA.jl}~\cite{besard2018juliagpu}, which extends the Julia compiler with GPU-specific constructs to generate GPU-compatible LLVM IR. While GPUs are traditionally used for compute-intensive workloads like graphics and deep learning, recent work has demonstrated their effectiveness for large batches of embarrassingly parallel problems, such as solving thousands of small differential equations in parallel~\cite{utkarsh2024automated}.

To support such workloads, \nlsolve{} includes \texttt{SimpleNonlinearSolve.jl}, a lightweight backend that provides highly efficient, non-allocating solvers for very small nonlinear systems. Algorithms like \texttt{Newton-Raphson} and \texttt{Trust-Region} are implemented using statically sized data structures via \texttt{StaticArrays.jl}, avoiding heap allocations and dynamic dispatch. This design makes them well-suited for GPU execution. These solvers integrate seamlessly with \texttt{KernelAbstractions.jl}~\cite{churavy2021juliagpu}, allowing users to embed them directly within GPU kernels and solve many independent nonlinear systems in parallel across GPU threads. In the example below, we solve the generalized Rosenbrock problem [\Cref{eq:gen_rosenbrock_function}] for $1024$ different initial conditions on CPU, AMD ROCm GPUs, and NVIDIA CUDA GPUs using identical code.
\begin{minted}[breaklines,escapeinside=||,mathescape=true, linenos, numbersep=3pt, gobble=2, frame=lines, fontsize=\ssmall, framesep=2mm]{julia}
using NonlinearSolve, StaticArrays
using KernelAbstractions, CUDA, AMDGPU

@kernel function parallel_nonlinearsolve_kernel!(result, @Const(prob), @Const(alg))
    i = @index(Global)
    prob_i = remake(prob; u0 = prob.u0[i])
    sol = solve(prob_i, alg)
    @inbounds result[i] = sol.u
    return nothing
end

function vectorized_solve(prob, alg; backend = CPU())
    result = KernelAbstractions.allocate(backend, eltype(prob.u0), length(prob.u0))

    groupsize = min(length(prob.u0), 1024)

    kernel! = parallel_nonlinearsolve_kernel!(backend, groupsize, length(prob.u0))
    kernel!(result, prob, alg)
    KernelAbstractions.synchronize(backend)

    return result
end

@generated function generalized_rosenbrock(x::SVector{N}, p) where {N}
    vals = ntuple(i -> gensym(string(i)), N)
    expr = []
    push!(expr, :($(vals[1]) = oneunit(x[1]) - x[1]))
    for i in 2:N
        push!(expr, :($(vals[i]) = 10.0 * (x[$i] - x[$i - 1] * x[$i - 1])))
    end
    push!(expr, :(@SVector [$(vals...)]))
    return Expr(:block, expr...)
end

u0 = @SVector [@SVector(rand(10)) for _ in 1:1024]
prob = NonlinearProblem(generalized_rosenbrock, u0)

vectorized_solve(prob, SimpleNewtonRaphson(); backend = CPU())
vectorized_solve(prob, SimpleNewtonRaphson(); backend = ROCBackend())
vectorized_solve(prob, SimpleNewtonRaphson(); backend = CUDABackend())
\end{minted}
The simpler solvers outperform the more general solvers in NonlinearSolve.jl significantly for small static problems [\Cref{fig:quadratic_problem_simplenonlinearsolve}]. Their high performance enables applications like massively parallel global optimization~\cite{utkarsh2024efficient} and parameter estimation problems, where solving many small independent nonlinear systems on the GPU is advantageous. \texttt{SimpleNonlinearSolve.jl} provides a portable, vendor-agnostic implementation that can target different GPU architectures like CUDA, ROCm, etc., with the same code.

\subsection{Automatic Sparsity Exploitation}\label{subsec:automatic-sparsity-detection}

Symbolic sparsity detection has a high overhead for smaller systems with well-defined sparsity patterns. We provide an approximate algorithm to determine the Jacobian sparsity pattern in those setups. We compute the dense Jacobian for $n$ randomly generated inputs to approximate the pattern. We take a union over the non-zero elements of Jacobian to obtain the sparsity pattern. As evident, computing the sparsity pattern costs $n$ times the cost of computing the dense Jacobian, typically via automatic forward mode differentiation.

\begin{figure}[t]
    \centering
    \adjustbox{trim=4em 2em 3em 3em, clip}{
        \includeinkscape[width=1.2\textwidth]{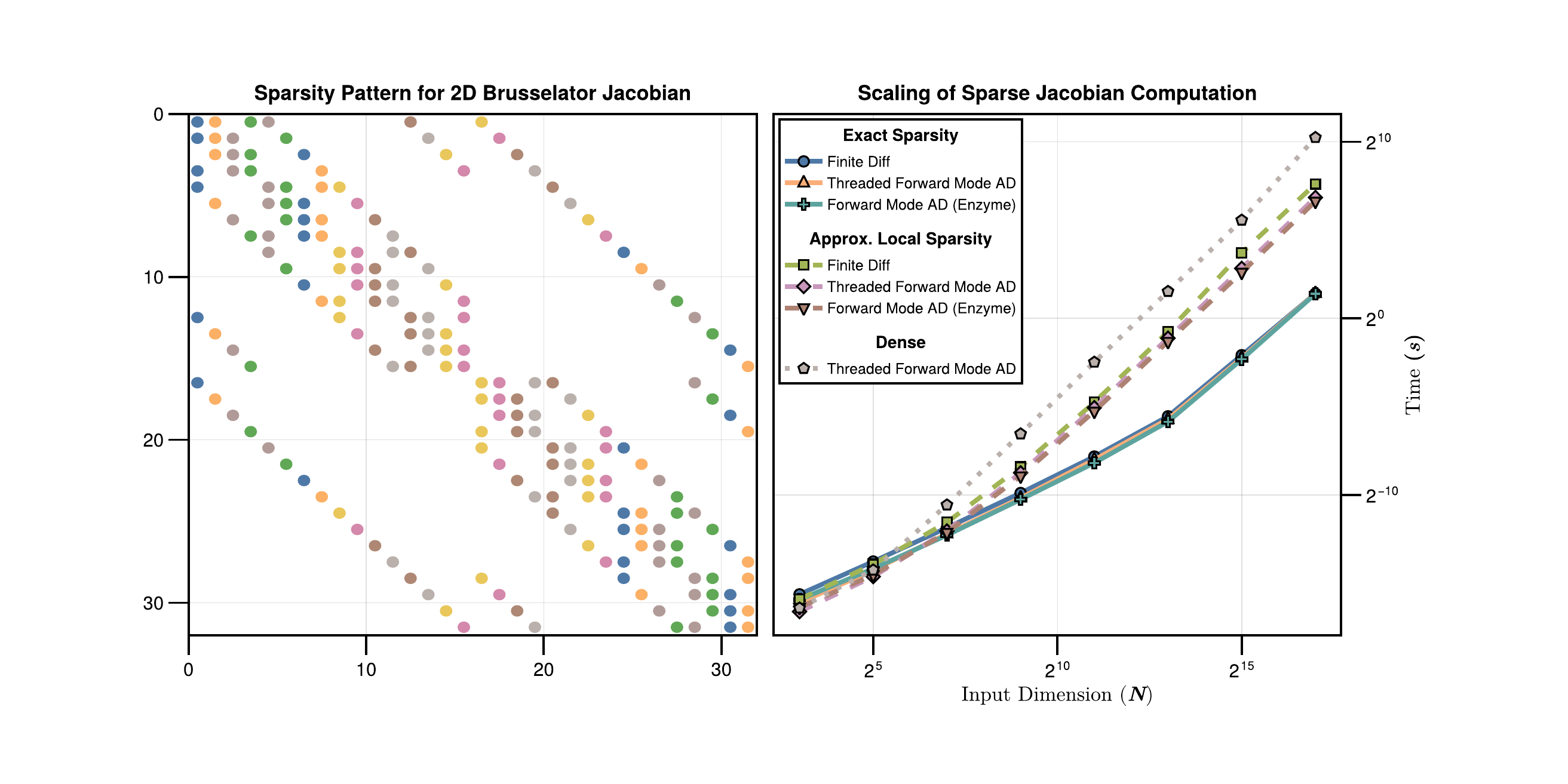}
    }
    \caption{\textbf{Automatic Sparsity Detection and Jacobian Computation for 2D Brusselator [\Cref{subsec:brusselator}]:} We benchmark the time to compute Jacobians using different modes of AD and FD under exact, approximate, and no sparsity assumptions. The left plot shows the exact sparsity pattern detected via matrix coloring, while the right plot compares the scaling behavior of different Jacobian computation methods as a function of input dimension. We include benchmarks using Enzyme-based AD~\cite{moses2020instead} to show that the choice between forward mode AD and finite differences has minimal impact on performance when sparsity is properly exploited. The dominant performance gains stem from exploiting sparsity structure, with exact sparsity consistently yielding better scaling for large systems.}\label{fig:brusselator_sparse_jacobian_scaling}
    \vspace{-1em}
\end{figure}

Approximate sparsity detection has poor scaling beyond a certain problem size, as evident from \Cref{fig:brusselator-scaling-benchmarks}. Similar to the shortcomings of other numerical sparsity detection software~\cite{giering2005generating, walther2009getting}, our method fails to accurately predict the sparsity pattern in the presence of state-dependent branches and might over-predict or under-predict sparsity due to floating point errors. Regardless, we observe in \Cref{fig:brusselator-scaling-benchmarks} that approximate sparsity detection is extremely efficient for moderately sized problems. In addition to computing the Jacobian faster, sparsity detection enables using sparse linear solvers that are significantly more efficient than solving the equivalent large dense linear systems [\Cref{subsec:brusselator}].

We rigorously benchmark the effects of sparsity detection algorithms in nonlinear root-finding in \Cref{subsec:brusselator}. Here, we will isolate the cost of sparsity detection on a 2D Brusselator with increasing discretization. In \Cref{fig:brusselator_sparse_jacobian_scaling}, we perform sparsity detection on the 2D Brusselator [\Cref{eq:bruss_pde}] and compute the Jacobian 10 times to resemble a realistic workload where detection is performed once and reused for later iterations. We can see that a threaded implementation of forward mode AD (\texttt{PolyesterForwardDiff}) can compute the dense Jacobian faster than most sparse Jacobian algorithms for small problems. After an initial threshold, using approximate sparsity detection techniques will outperform other techniques. Finally, for large systems, using exact symbolic sparsity detection followed by colored AD is the most efficient.

\subsection{Generalized Jacobian-Free Nonlinear Solvers using Krylov Methods}\label{subsec:krylov-methods}

\nlsolve{} supports Jacobian-Free solvers such as Newton-Krylov and Dogleg with Krylov methods via the iterative linear solvers in \texttt{LinearSolve.jl}. A trait-based interface detects whether the selected linear solver supports matrix-free computation. If supported, we construct a \texttt{JacobianOperator} that evaluates $\mathcal{J}x$ via forward-mode AD and $x^T \mathcal{J}$ via reverse-mode AD. For example, \texttt{NewtonRaphson} and \texttt{TrustRegion} --- the two most commonly used algorithms in \nlsolve{} --- compute descent directions as follows:
\begin{align}
    \mathcal{J}_k \delta u_{\text{newton}} = -f(u_k, \theta)\label{eq:newton_descent_direction}
\end{align}
\Cref{eq:newton_descent_direction} is solved using iterative methods, with $\mathcal{J}_k \delta u_{\text{newton}}$ evaluated via forward-mode AD. For Dogleg, both the Newton and steepest descent directions are required. The steepest descent direction is:
\begin{align}
    \delta u_{\text{sd}} = -\mathcal{J}_k^T f(u_k, \theta)
\end{align}
This is computed directly using reverse-mode AD without materializing $\mathcal{J}_k$. The \texttt{JacobianOperator} supports both $\mathcal{J}x$ and $x^T \mathcal{J}$, making it compatible with least-squares Krylov solvers like \texttt{LSMR}~\cite{fong2010lsmr}. For some Krylov methods, convergence depends on preconditioning, which may require a concrete Jacobian. To support this, we provide the \texttt{concrete\_jac} option, allowing users to force materialization of $\mathcal{J}$ when necessary. In this mode, both Jacobian-vector and vector-Jacobian products are computed explicitly. In \Cref{subsec:brusselator}, we demonstrate Jacobian-Free Newton and Dogleg methods using \texttt{GMRES}~\cite{saad1986gmres}, preconditioned with \texttt{IncompleteLU.jl}~\cite{Stoppels2024haampie} and \texttt{AlgebraicMultigrid.jl}~\cite{Anantharaman2025JuliaLinearAlgebra}. For large-scale systems, Krylov methods significantly outperform direct methods [\Cref{fig:brusselator-scaling-benchmarks-krylovmethods}] and benefit from our sparse Jacobian infrastructure, which enables efficient preconditioner construction.

\begin{table}[t]
    \centering
    \begin{tabular}{l l l l}
        \toprule
        \makecell{Framework} & \makecell{Version} & \makecell{Backend Language} & \makecell{Build Script} \\
        \midrule
        \texttt{NonlinearSolve}    & $4.6.0$  & \texttt{Julia} & N/A \\
        \texttt{NLsolve}           & $4.5.1$  & \texttt{Julia} & N/A \\
	    \texttt{Sundials}          & $5.2$    & \texttt{C}     & \href{https://github.com/JuliaPackaging/Yggdrasil/blob/bb48b83737ee7d4b6c37efe1631f6232ce664b1b/S/Sundials/Sundials\%405/build\_tarballs.jl}{build\_tarballs.jl} \\
        \texttt{CMINPACK}          & $1.3.11$ & \texttt{C}     & \href{https://github.com/JuliaPackaging/Yggdrasil/blob/bb48b83737ee7d4b6c37efe1631f6232ce664b1b/C/cminpack/build\_tarballs.jl}{build\_tarballs.jl} \\
        \texttt{PETSc}             & $3.22.0$ & \texttt{C}     & \href{https://github.com/JuliaPackaging/Yggdrasil/blob/bb48b83737ee7d4b6c37efe1631f6232ce664b1b/P/PETSc/build\_tarballs.jl}{build\_tarballs.jl} \\
        \bottomrule
    \end{tabular}
    \vspace{1em}
    \caption{Versions for the frameworks used in the experiments and the build scripts for compiling shared libraries for non-Julia libraries.}\label{tab:library_versions}
    \vspace{-2em}
\end{table}

\begin{table}[t]
    \begin{adjustbox}{width=0.9\textwidth,center}
    \begin{tabular}{l l l l l l c}
        \toprule
        \makecell{Method} & \makecell{Line Search} & \makecell{Trust Region\\Scheme} & \makecell{Linear\\Solver} & \makecell{Framework} & \makecell{AD / FD} & \makecell{Number of\\Problems Solved} \\
        \midrule
        Newton Raphson & - & - & LU & \nlsolve{} & AD & $23$ \\
        Newton Raphson & Hager Zhang & - & LU & \nlsolve{} & AD & $20$ \\
        Newton Raphson & More Thuente & - & LU & \nlsolve{} & AD & $22$ \\
        Newton Raphson & Back Tracking & - & LU & \nlsolve{} & AD & $22$ \\
        Trust Region   & - & Simple & LU & \nlsolve{} & AD & $21$ \\
        Trust Region   & - & Hei & LU & \nlsolve{} & AD & $21$ \\
        Trust Region   & - & Fan & LU & \nlsolve{} & AD & $21$ \\
        Trust Region   & - & NLsolve & LU & \nlsolve{} & AD & $21$ \\
        Trust Region   & - & Bastin & LU & \nlsolve{} & AD & $22$ \\
        Trust Region   & - & NocedalWright & LU & \nlsolve{} & AD & $21$ \\
        Trust Region   & - & Yuan & LU & \nlsolve{} & AD & $20$ \\
        Levenberg Marquardt   & - & - & QR & \nlsolve{} & AD & $21$ \\
        Levenberg Marquardt   & - & - & Cholesky & \nlsolve{} & AD & $21$ \\
        L.M. (No Geodesic Accln.)  & - & - & QR & \nlsolve{} & AD & $23$ \\
        L.M. (No Geodesic Accln.)  & - & - & Cholesky & \nlsolve{} & AD & $23$ \\
        Pseudo Transient & - & - & LU & \nlsolve{} & AD & $20$ \\
        Newton Raphson & - & - & LU & \texttt{Sundials} & FD & $22$ \\
        Newton Raphson & Back Tracking & - & LU & \texttt{Sundials} & FD & $22$ \\
        Newton Raphson & - & - & LU & \texttt{NLsolve.jl} & AD & $22$ \\
        Trust Region   & - & NLsolve & LU & \texttt{NLsolve.jl} & AD & $22$ \\
        Modified Powell & - & - & QR & \texttt{MINPACK} & FD & $18$ \\
        Levenberg Marquardt & - & - & QR & \texttt{MINPACK} & FD & $17$ \\
        \bottomrule
    \end{tabular}
    \end{adjustbox}
    \vspace{1em}
    \caption{\textbf{Solver Performance on 23 Small-Scale Nonlinear Root-Finding Benchmarks.} We evaluate various solver configurations from \nlsolve{}, \texttt{Sundials}, \texttt{CMINPACK}, and \texttt{NLsolve.jl} on a suite of 23 nonlinear problems. \texttt{PETSc} is excluded from these benchmarks due to its relatively high overhead for small-scale problems (2–10 variables), as demonstrated in \Cref{subsec:brusselator} and \Cref{subsec:battery_model}.}\label{tab:23_test_problems}
    \vspace{-2em}
\end{table}

\section{Results}\label{sec:results}

We evaluate our solvers on three numerical experiments and benchmark them against other nonlinear equation solvers -- \texttt{PETSc SNES}~\cite{balay2019petsc}, \texttt{NLsolve.jl}~\cite{Optim.jl-2018}, \texttt{Sundials}~\cite{gardner2022sundials, hindmarsh2005sundials}, and \texttt{CMINPACK}~\cite{cminpack} (a modern C/C++ rewrite of \texttt{MINPACK}). We run all our benchmarks on \texttt{Julia 1.10} on a single node with \texttt{128 × AMD EPYC 7502 32-Core Processor} with \texttt{128 Julia Threads}. All solvers use the $L_\infty$-norm of the residual to determine convergence and run until convergence or a maximum of $1000$ iterations. \Cref{tab:library_versions} summarizes the versions of the frameworks used in the experiments. The code for our numerical experiments is publicly available on \href{https://github.com/SciML/SciMLBenchmarks.jl/tree/2471352ed69d613424d76e7084faf2915227f4a7/benchmarks/NonlinearProblem}{GitHub}.

\subsection{Evaluation on 23 Small-Scale Nonlinear Root-Finding Benchmarks}\label{subsec:23_test_problems}

The \href{https://web.archive.org/web/20240326025519/https://people.sc.fsu.edu/~jburkardt/f_src/test_nonlin/test_nonlin.html}{23 Test Problems} suite is a standard collection of small-scale nonlinear systems, each with $2$ to $10$ variables. These problems are commonly used to evaluate the performance, robustness, and accuracy of nonlinear solvers. The original FORTRAN90 implementation includes routines for initializing each problem, computing function values, and evaluating the Jacobian. We provide a Julia-based version of this suite via \href{https://github.com/SciML/DiffEqProblemLibrary.jl/blob/69655a67a1a0d8575ec87e05a2013d5a3a27628e/lib/NonlinearProblemLibrary}{\texttt{NonlinearProblemLibrary.jl}}, where Jacobians are computed automatically using AD or finite differencing.

\begin{figure}[t]
    \centering
    \adjustbox{trim=4em 3.5em 3em 6em, clip}{
        \includeinkscape[width=1.15\textwidth]{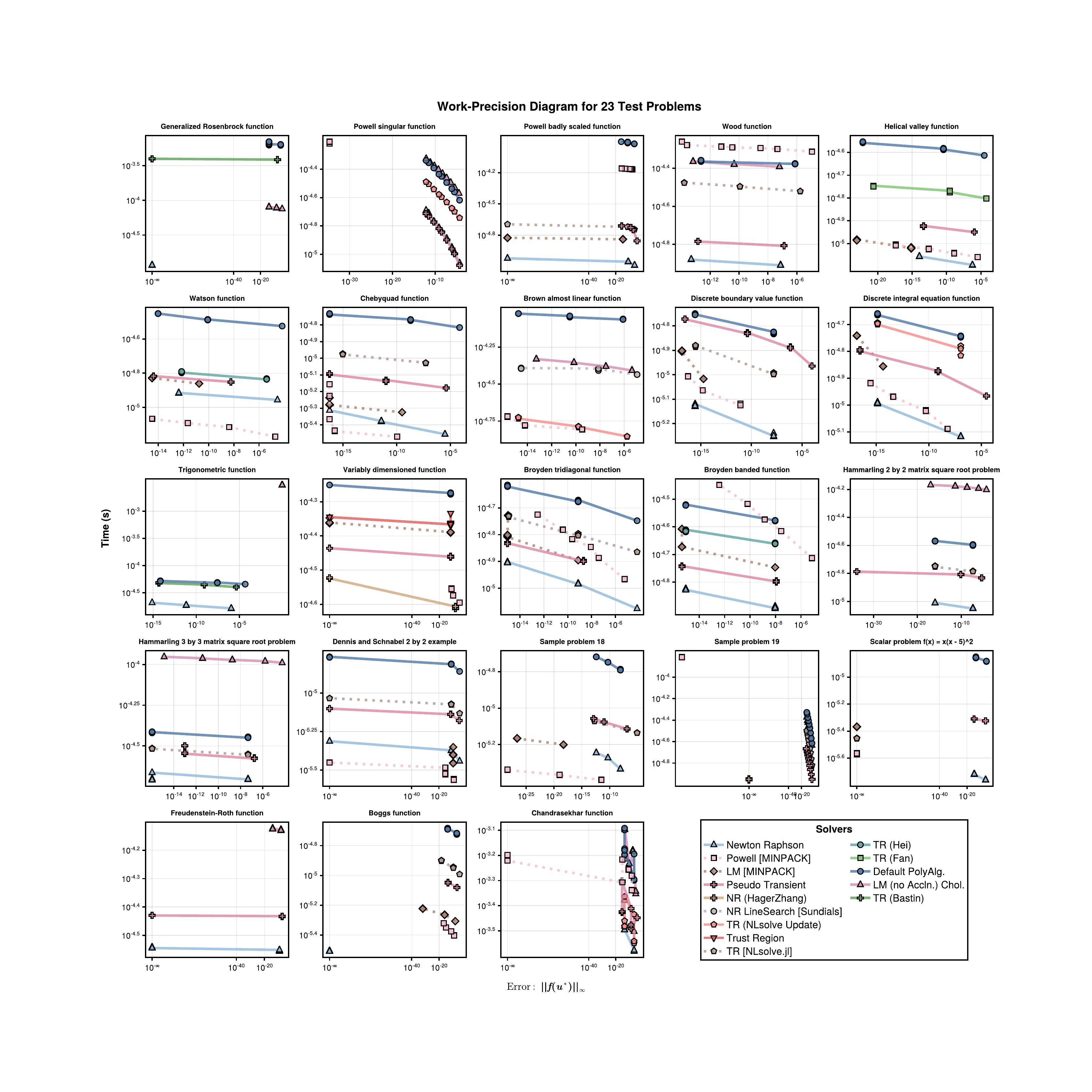}
    }
    \vspace{-3em}
    \caption{\textbf{Work-Precision Diagrams for 23 Small-Scale Nonlinear Root-Finding Benchmarks.} Each subplot shows solver runtime versus residual error $\|f(u^*)\|$ for one of the 23 test problems. While solver performance varies by problem and algorithm, most configurations perform comparably across the suite. \texttt{NonlinearSolve.jl} successfully solves all problems and demonstrates strong performance overall. These plots illustrate the trade-offs between accuracy and execution time across solver strategies, with the \texttt{Newton-Raphson} implementation in \nlsolve{} generally outperforming most other frameworks. For clarity, each subplot includes only the fastest-performing implementation within each solver category. For instance, if \texttt{Newton-Raphson} from \texttt{NonlinearSolve.jl} achieves the lowest median runtime among all Newton-based methods (including those with line search), only that configuration is shown. Full results and detailed plots for each problem are available on our \href{https://docs.sciml.ai/SciMLBenchmarksOutput/dev/NonlinearProblem/nonlinear_solver_23_tests/}{benchmark results page}.}
    \label{fig:23test_problems_summary}
    \vspace{-1em}
\end{figure}

\Cref{tab:23_test_problems} lists all the solver configurations included in this experiment, and \Cref{fig:23test_problems_summary} summarizes their performance. To ensure a fair comparison and control for differences in linear solver efficiency, all methods are configured to use \texttt{LAPACK}-based linear solvers. While many solvers perform similarly across a range of problems, \texttt{NonlinearSolve.jl} demonstrates greater flexibility through its poly-algorithmic approach and broader selection of solver strategies, successfully solving all 23 test cases. Other packages, including \texttt{Sundials}, \texttt{MINPACK}, and \texttt{NLsolve.jl}, also perform well on a majority of the problems and serve as robust alternatives, although they exhibit slightly higher overhead on small-scale systems in this benchmark.

\subsection{Initializing the Doyle-Fuller-Newman (DFN) Battery Model}\label{subsec:battery_model}

Initializing DAEs is a critical step in the numerical solution process. It ensures consistency and well-posedness of the problem. Poor initialization can lead to solver divergence, inaccurate results, or total failure. Initialization involves computing consistent initial conditions for both differential and algebraic variables such that the algebraic constraints are satisfied at the initial time. This guarantees that the numerical solver starts from a physically meaningful state without violating system constraints.

We consider a 32-dimensional DAE initialization problem based on the \href{https://help.juliahub.com/batteries/stable/}{Doyle-Fuller-Newman (DFN) battery model}~\cite{doyle1993modeling, rackauckas2022composing}. The initial state is derived from the open-circuit voltage (OCV) profile under high current charge conditions~\cite{berliner2021methods, rackauckas2022composing}. This initialization is challenging due to the stiffness and nonlinearity of the coupled electrochemical equations, where even small deviations in guessed states can lead to solver failure, particularly because the Jacobian becomes ill-conditioned (with a condition number of $\approx 10^{14}$), difficult to invert under high current rates~\cite{berliner2021methods}.

\begin{table}[t]
    \begin{adjustbox}{width=0.8\textwidth,center}
    \begin{tabular}{l l l l l l c}
        \toprule
        \makecell{Method} & \makecell{Line Search} & \makecell{Trust Region\\Scheme} & \makecell{Linear\\Solver} & \makecell{Framework} & \makecell{AD / FD} & \makecell{Successful} \\
        \midrule
        Newton Raphson & - & - & LU & \nlsolve{} & AD & \cmark \\
        Newton Raphson & Hager Zhang & - & LU & \nlsolve{} & AD & \xmark \\
        Newton Raphson & More Thuente & - & LU & \nlsolve{} & AD & \xmark \\
        Newton Raphson & Back Tracking & - & LU & \nlsolve{} & AD & \xmark \\
        Newton Krylov  & - & - & GMRES & \nlsolve{} & AD & \xmark \\
        Trust Region   & - & Simple & LU & \nlsolve{} & AD & \cmark \\
        Trust Region   & - & Hei & LU & \nlsolve{} & AD & \cmark \\
        Trust Region   & - & Fan & LU & \nlsolve{} & AD & \cmark \\
        Trust Region   & - & NLsolve & LU & \nlsolve{} & AD & \cmark \\
        Trust Region   & - & Bastin & LU & \nlsolve{} & AD & \cmark \\
        Trust Region   & - & NocedalWright & LU & \nlsolve{} & AD & \xmark \\
        Trust Region   & - & Yuan & LU & \nlsolve{} & AD & \xmark \\
        Levenberg Marquardt   & - & - & QR & \nlsolve{} & AD & \cmark \\
        Levenberg Marquardt   & - & - & Cholesky & \nlsolve{} & AD & \xmark \\
        L.M. (No Geodesic Accln.)  & - & - & QR & \nlsolve{} & AD & \cmark \\
        L.M. (No Geodesic Accln.)  & - & - & Cholesky & \nlsolve{} & AD & \xmark \\
        Newton Raphson & - & - & LU & \texttt{Sundials} & FD & \cmark \\
        Newton Krylov  & - & - & GMRES & \texttt{Sundials} & FD & \xmark \\
        Newton Raphson & Back Tracking & - & LU & \texttt{Sundials} & FD & \xmark \\
        Newton Raphson & - & - & LU & \texttt{NLsolve.jl} & AD & \cmark \\
        Trust Region   & - & NLsolve & LU & \texttt{NLsolve.jl} & AD & \cmark \\
        Modified Powell & - & - & QR & \texttt{MINPACK} & FD & \xmark \\
        Levenberg Marquardt & - & - & QR & \texttt{MINPACK} & FD & \xmark \\
        Newton Raphson & - & - & LU & \texttt{PETSc} & FD & \xmark \\
        Newton Raphson & - & - & QR & \texttt{PETSc} & FD & \xmark \\
        Newton Raphson & Back Tracking & - & LU & \texttt{PETSc} & FD & \xmark \\
        Newton Raphson & Back Tracking & - & QR & \texttt{PETSc} & FD & \xmark \\
        Trust Region & - & - & LU & \texttt{PETSc} & FD & \xmark \\
        Newton Krylov & - & - & GMRES & \texttt{PETSc} & FD & \xmark \\
        \bottomrule
    \end{tabular}
    \end{adjustbox}
    \vspace{1em}
    \caption{\textbf{Initializing the DFN Battery DAE Model at a High Current}: Algorithms and configurations evaluated on this problem. We exclude Quasi-Newton methods, as they are known to fail~\cite{berliner2021methods}. A solution is marked successful if the $L_\infty$-norm between the computed and reference solutions is less than $1\times 10^{-4}$.}\label{tab:battery_problem_results}
    \vspace{-2em}
\end{table}

\begin{figure}[bthp]
    \centering
    \adjustbox{trim=4em 3em 3em 5em, clip}{
        \includeinkscape[width=1.15\textwidth]{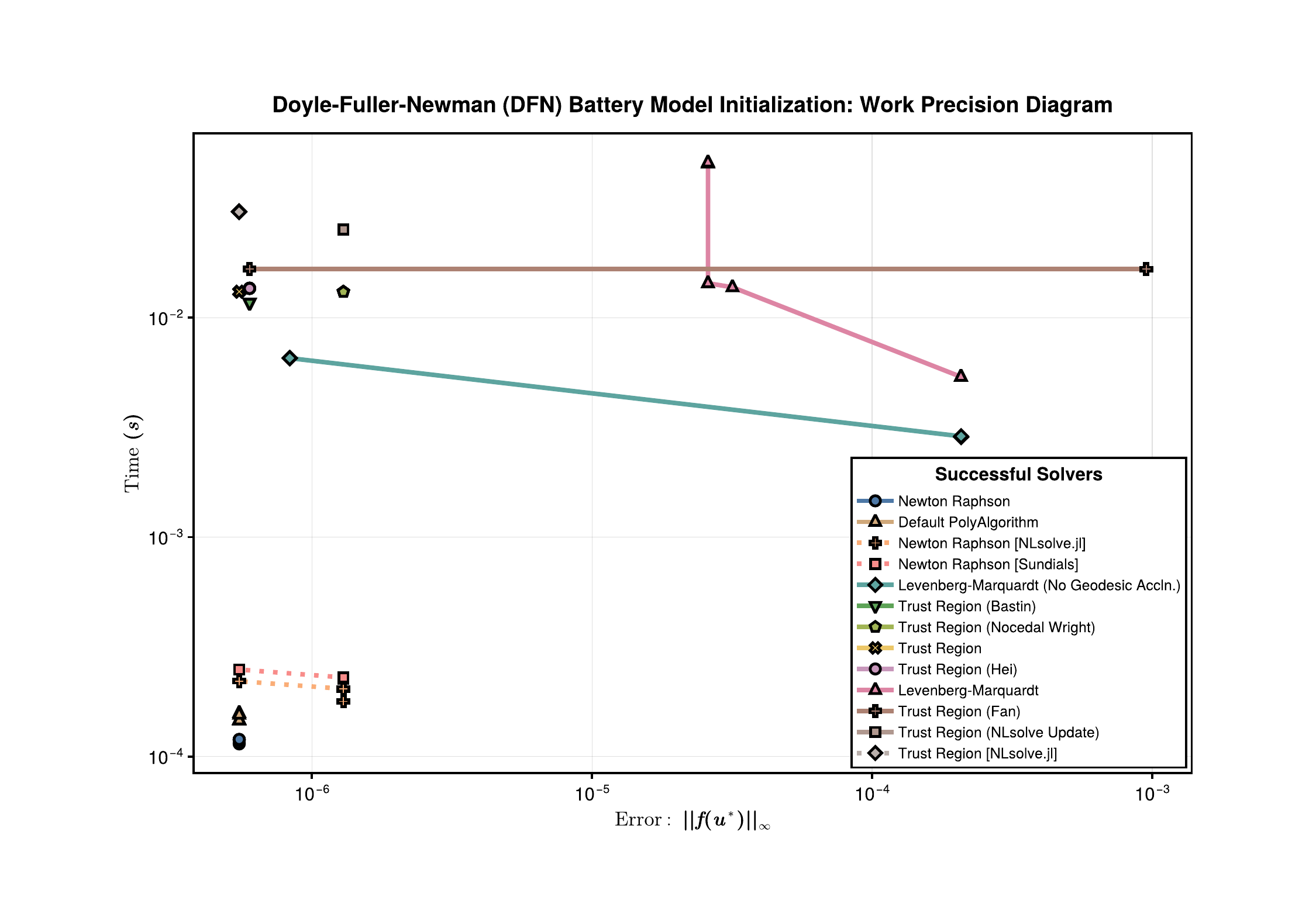}
    }
    \caption{\textbf{Work-Precision Diagram for DFN Battery DAE Initialization at High Current}: Among all frameworks tested, \nlsolve{} solvers --- \texttt{TrustRegion}, \texttt{NewtonRaphson}, and \texttt{LevenbergMarquardt} --- consistently solve the problem and achieve the best work-precision tradeoff. These outperform the same solvers implemented in \texttt{NLsolve.jl} \& \texttt{Sundials KINSOL}. All tested \texttt{MINPACK} and \texttt{PETSc} solvers fail to converge below all tested target tolerances.}\label{fig:battery-model-benchmarks}
\end{figure}

\Cref{fig:battery-model-benchmarks} illustrates the work-precision tradeoff for all tested solvers. \nlsolve{} methods solve the problem reliably, with \texttt{NewtonRaphson} providing the most efficient convergence. In contrast, solvers from \texttt{MINPACK} and \texttt{PETSc} fail across the board, underscoring the difficulty of this DAE system. %

\begin{figure}[t]
    \centering
    \adjustbox{trim=4em 3em 3em 5em, clip}{
        \includeinkscape[width=1.15\textwidth]{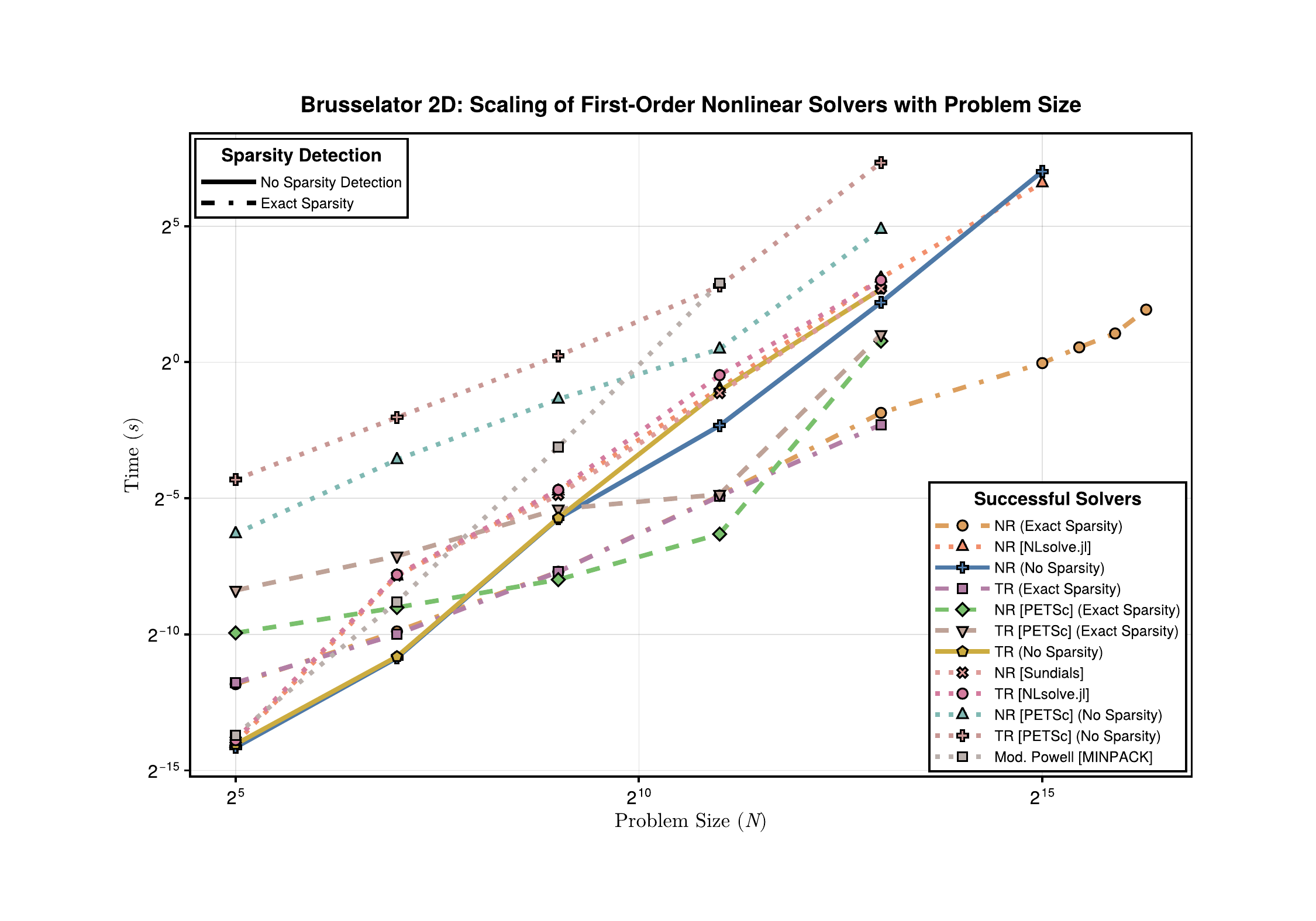}
    }
    \caption{\textbf{Brusselator 2D Scaling of First-Order Nonlinear Solvers.} This plot compares the runtime scaling of first-order solvers across different frameworks for increasing problem sizes. \texttt{NonlinearSolve.jl}'s \texttt{Newton-Raphson} and \texttt{Trust-Region} implementations show strong problem-size scaling performance, with trust-region methods significantly outperforming the widely used Modified Powell solver (a trust-region based method) from \texttt{MINPACK}, commonly used in tools like \texttt{SciPy} and \texttt{MATLAB}. Solvers in \texttt{PETSc} that use sparsity patterns detected by \texttt{NonlinearSolve.jl} are competitive for medium-scale problems but fail to converge below a target tolerance of $10^{-6}$ within $10$-minutes for larger sizes. Overall, \texttt{NonlinearSolve.jl} demonstrates robust performance and improved scalability, particularly when leveraging sparsity detection.}\label{fig:brusselator-scaling-benchmarks}
\end{figure}

\begin{figure}[t]
    \centering
    \adjustbox{trim=4em 3em 3em 5em, clip}{
        \includeinkscape[width=1.15\textwidth]{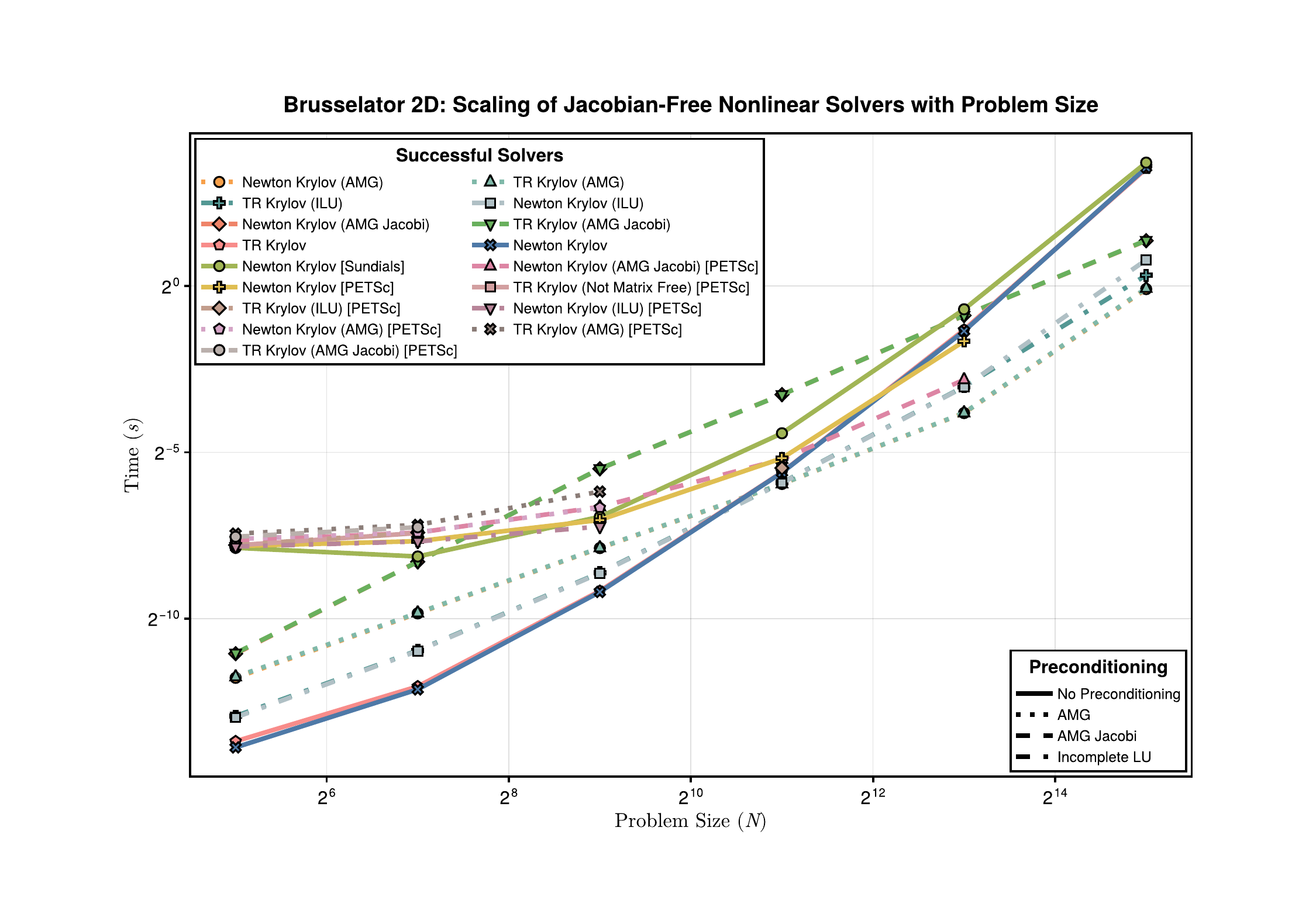}
    }
    \caption{\textbf{Brusselator 2D Problem-Size Scaling with Jacobian-Free Krylov Solvers.} \texttt{NonlinearSolve.jl} integrates with \texttt{LinearSolve.jl} to support Jacobian-free Krylov methods such as \texttt{GMRES}, \texttt{CG}, and \texttt{LSMR}, specified via \texttt{linsolve = KrylovJL\_$\ast$}. These solvers avoid explicit Jacobian construction, reducing overhead for large systems as seen in \Cref{fig:brusselator-scaling-benchmarks}. Preconditioning is optionally applied to improve convergence, and consistent Krylov solver parameters (no restarts, fixed Krylov dimension) are used across frameworks for fairness. Only configurations that achieve a residual below $10^{-4}$ are shown. \texttt{NonlinearSolve.jl} demonstrates faster runtimes across both preconditioned and no-preconditioned cases.}\label{fig:brusselator-scaling-benchmarks-krylovmethods}
\end{figure}

\begin{figure}[tbhp]
    \centering
    \adjustbox{trim=4em 3em 3em 5em, clip}{
        \includeinkscape[width=1.15\textwidth]{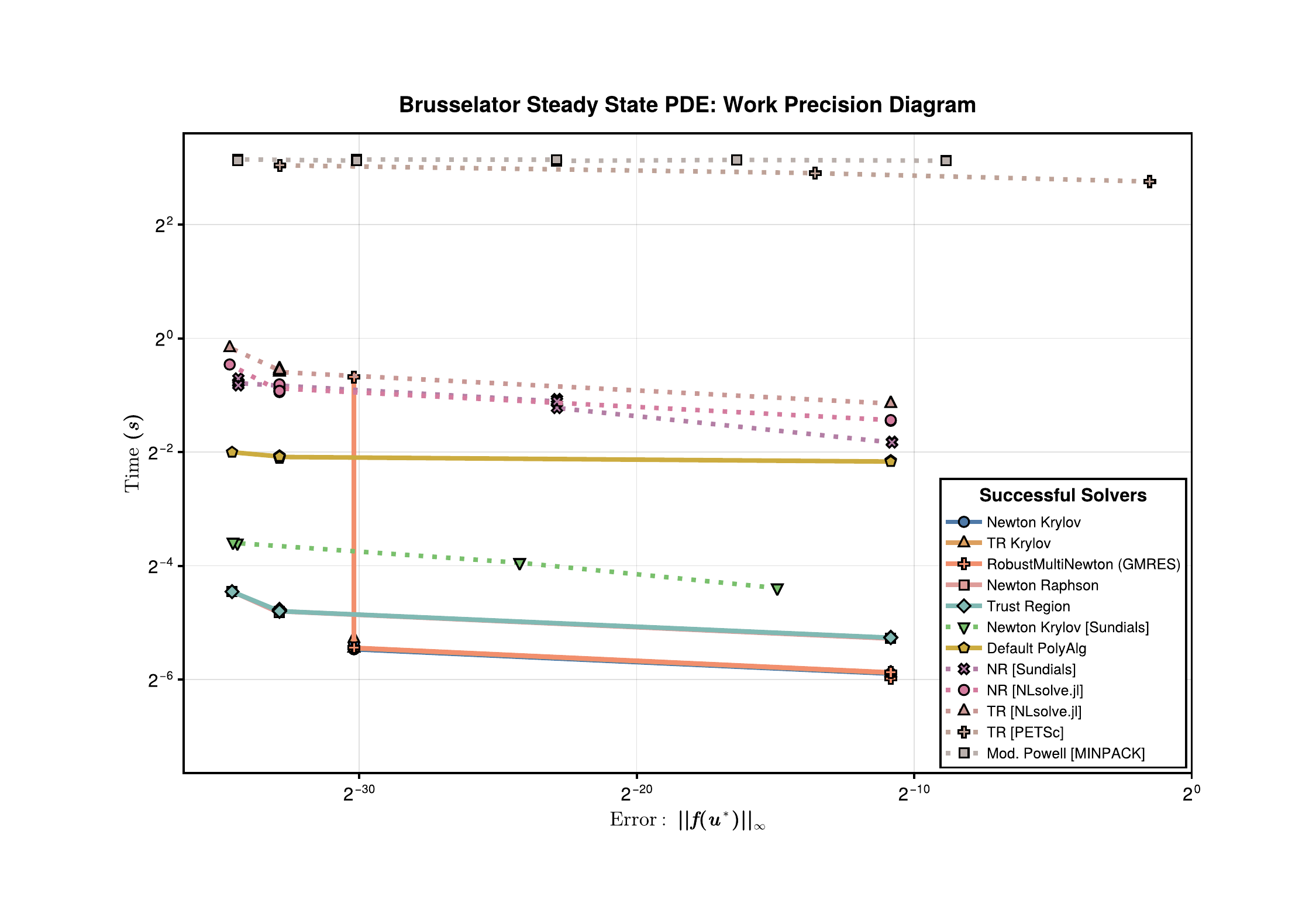}
    }
    \caption{\textbf{Work-Precision Diagram for the Brusselator Steady-State PDE.} This plot compares solver runtime against residual error $\|f(u^*)\|_\infty$ for computing the steady state of the 2D Brusselator system. \texttt{NonlinearSolve.jl} solvers, including Jacobian-free Krylov and poly-algorithmic variants like \texttt{RobustMultiNewton}, generally exhibit faster runtimes compared to their counterparts in \texttt{MINPACK}, \texttt{Sundials}, \texttt{PETSc}, and \texttt{NLsolve.jl}. Among external libraries, \texttt{Sundials KINSOL} with GMRES (without Jacobian reuse) achieves the best performance.}\label{fig:brusselator-wpd}
\end{figure}

\subsection{Large Ill-Conditioned Nonlinear Brusselator System}\label{subsec:brusselator}

Solving ill-conditioned nonlinear systems requires specializing the linear solver on properties of the Jacobian to cut down on the $\mathcal{O}(n^3)$ linear solve and the $\mathcal{O}(n^2)$ back-solves. We will use the sparsity detection and Jacobian-free solvers from \nlsolve{} to solve the steady state stiff Brusselator PDE \cite[Section 4.2]{wanner1996solving}. The Brusselator PDE with  is defined as follows:
\begin{equation}
    \begin{aligned}
        0 &= 1 + u^2 v - 4.4u + \alpha \nabla^2 u + f(x, y, t), \\
        0 &= 3.4u - u^2 v + \alpha \nabla^2 v.
    \end{aligned}\label{eq:bruss_pde}
\end{equation}
where
\begin{align}
    f(x, y, t) &= \begin{cases}
        5 & \quad \text{if } (x-0.3)^2+(y-0.6)^2 \leq 0.1^2 \text{ and } t \geq 1.1 \\
        0 & \quad \text{else}
    \end{cases} \\
    \alpha &= 10 (N - 1)^2 \quad \text{ N is the discretization}
\end{align}
and the initial conditions are
\begin{align}
    u(x, y, 0) & = 22\cdot (y(1-y))^{3/2} \\
    v(x, y, 0) & = 27\cdot (x(1-x))^{3/2}
\end{align}
with the periodic boundary condition
\begin{align}
    u(x+1,y,t) & = u(x,y,t) \\
    u(x,y+1,t) & = u(x,y,t)
\end{align}
We analyze our solvers for this problem in 3 parts. Firstly, \Cref{fig:brusselator-scaling-benchmarks} compares the scaling of \texttt{NewtonRaphson} and \texttt{TrustRegion} between frameworks. \nlsolve{} leverages sparsity detection and colored sparse matrix algorithms to accelerate Jacobian construction and is consistently faster than all existing solvers here. Most of these frameworks have built-in sparsity support and, unsurprisingly, don't scale well with increasing problem size. However, with increasing discretization, explicitly forming Jacobians eventually becomes infeasible. \nlsolve{} allows automatically building sparse Jacobians, that can be consumed by frameworks like \texttt{PETSc} as shown in the plot. \Cref{fig:brusselator-scaling-benchmarks-krylovmethods} demonstrates \nlsolve{} Jacobian-free Krylov methods with preconditioning that solve the linear system without constructing the Jacobian. Newton-Krylov methods from \texttt{Sundials} and \texttt{PETSc} are competitive with \nlsolve{} at larger problem sizes. \texttt{MINPACK} and \texttt{NLsolve.jl} don't support Jacobian-Free methods.  Finally, in \Cref{fig:brusselator-wpd}, we present a Work-Precision Diagram for the discretized Brusselator on a $32 \times 32$ grid, showing that our methods are an order of magnitude faster than the ones in other software.

\section{Conclusion}
\label{sec:conclusion}

Solving systems of nonlinear equations is a fundamental challenge that arises across many scientific domains. This paper presented \nlsolve, a high-performance and robust open-source solver for nonlinear systems implemented natively in the Julia programming language. Through extensive numerical experiments on benchmark problems, applications from practical DAE initialization problems, and problem size scaling tests, we have demonstrated the advantages of using \nlsolve{} compared to existing software tools.

Key strengths of \nlsolve{} include its flexible unified API for rapidly experimenting with different solver options, smart automatic algorithm selection for balancing speed and robustness, specialized non-allocating kernels for small systems, automatic sparsity exploitation, and support for Jacobian-free Krylov methods. These features enable \nlsolve{} to solve challenging nonlinear problems, including cases where standard solvers fail while attaining high performance.

\begin{acks}

This material is based upon work supported by the Department of Energy, National Nuclear Security Administration under Award Number DE-NA0003965.This report was prepared as an account of work sponsored by an agency of the United States Government. Neither the United States Government nor any agency thereof, nor any of their employees, makes any warranty, express or implied, or assumes any legal liability or responsibility for the accuracy, completeness, or usefulness of any information, apparatus, product, or process disclosed, or represents that its use would not infringe privately owned rights. Reference herein to any specific commercial product, process, or service by trade name, trademark, manufacturer, or otherwise does not necessarily constitute or imply its endorsement, recommendation, or favoring by the United States Government or any agency thereof. The views and opinions of authors expressed herein do not necessarily state or reflect those of the United States Government or any agency thereof. This material is based upon work supported by the National Science Foundation under grant no. OAC-2103804 , no. OSI-2029670, no. DMS-2325184, no. PHY-2028125. This material was supported by The Research Council of Norway and Equinor ASA through the Research Council project ”308817 - Digital wells for optimal production and drainage”. Research was sponsored by the United States Air Force Research Laboratory and the United States Air Force Artificial Intelligence Accelerator and was accomplished under Cooperative Agreement Number FA8750-19-2-1000. The views and conclusions contained in this document are those of the authors and should not be interpreted as representing the official policies, either expressed or implied, of the United States Air Force or the U.S. Government. The U.S. Government is authorized to reproduce and distribute reprints for Government purposes, notwithstanding any copyright notation herein. The authors would like to thank DARPA for funding this work through the Automating Scientific Knowledge Ex-traction and Modeling (ASKEM) program, Agreement No. HR0011262087. The views, opinions, and/or findings expressed are those of the authors and should not be interpreted as representing the official views or policies of the Department of Defense or the U.S. Government.

\end{acks}

\bibliographystyle{ACM-Reference-Format}
\bibliography{references}

\end{document}